%% file: main.tex
\documentclass[preprint,12pt]{elsarticle}

\usepackage{amssymb}
 \usepackage{amsthm}

\input{mycustom}

\journal{Computer Methods in Applied Mechanics and Engineering}

\begin{document}

\begin{frontmatter}



\title{Robust parallel nonlinear solvers for implicit time discretizations  of the Bidomain equations}


\author[santiago]{Nicol{\'a}s A. Barnafi}

\affiliation[santiago]{organization={Centro de Modelamiento Matemático},
            addressline={Av. Beauchef 851}, 
            city={Santiago},
            postcode={8370456}, 
            state={},
            country={Chile}}
        
\author[pavia]{Ngoc Mai Monica Huynh}
\author[pavia]{Luca F. Pavarino}

\affiliation[pavia]{organization={Department of Mathematics, University of Pavia},
	addressline={via Ferrata 1}, 
	city={Pavia},
	postcode={27100}, 
	state={},
	country={Italy}}

\author[milano]{Simone Scacchi}

\affiliation[milano]{organization={Department of Mathematics, University of Milano},
	addressline={via Saldini 50}, 
	city={Milano},
	postcode={20133}, 
	state={},
	country={Italy}}

\begin{abstract}
In this work, we study the convergence and performance of nonlinear solvers for the Bidomain equations after decoupling the ordinary and partial differential equations of the cardiac system. Firstly, we provide a rigorous proof of the global convergence of  Quasi-Newton methods, such as BFGS, and nonlinear Conjugate-Gradient methods, such as Fletcher--Reeves, for the Bidomain system, by analyzing an auxiliary variational problem under physically reasonable hypotheses. Secondly, we compare several nonlinear Bidomain solvers in terms of execution time, robustness with respect to the data and parallel scalability. Our findings indicate that Quasi-Newton methods are the best choice for nonlinear Bidomain systems, since they exhibit faster convergence rates compared to standard Newton-Krylov methods, while maintaining robustness and scalability. Furthermore, first-order methods also demonstrate competitiveness and serve as a viable alternative, particularly for matrix-free implementations that are well-suited for GPU computing.

\end{abstract}




\end{frontmatter}


\section{Introduction} 

The complexity and significance of the human heart's functioning have long captivated researchers. 
Mathematical modeling plays an important role in understanding the underlying mechanisms governing the heart physiological and pathological conditions, as well as in the development of new tools, such as digital twins \cite{bjornsson2020digital}.
Many cardiac functions can be described mathematically by means of systems of partial differential equations (PDEs) and ordinary differential equations (ODEs), which are coupled together in order to describe different biological events, e.g. electrophysiology \cite{franzone2006computational}, muscle contraction and relaxation \cite{smith2004multiscale} and blood circulation \cite{dede2019computational, di2021computational, piersanti20213d}. It has also been extended to an optimal control framework to estimate activation sites \cite{nagaiah2011numerical,kunisch2013optimal}.

The efficient numerical simulation of these phenomena requires a balance between accuracy of the solution and computational efficiency.
Numerical methods needed to describe and computationally reproduce the many interactions between macroscopic and microscopic events often yield large nonlinear algebraic systems of equations, with millions of degrees of freedom. To this end, the development of robust, efficient, and scalable parallel solvers is of great importance. Many studies have focused on the development of parallel solvers for electro-mechanical models \cite{colli2018numerical} or on the coupling of the different physics involved \cite{quarteroni2017integrated}, as well as the theoretical and numerical analysis of the considered model \cite{pennacchio2005multiscale, barnafi2021mathematical}.

In this work, we focus on the Bidomain equations, a system of parabolic PDEs \cite{franzone2006computational} which models the propagation of the electric signal in the cardiac tissue, known as myocardium. This system is usually coupled with a nonlinear reaction term to a system of ODEs which, in turn, represents the ionic dynamics at a cellular level.
This system has been extensively studied within the framework of parallel solvers and preconditioners. For instance, by addressing a staggered decoupled time discretization, where at each time step the ODEs are solved before the PDEs, many studies have been devoted to the development and numerical validation of efficient preconditioners, such as Refs. \cite{chen2019splitting, munteanu2009sisc, pennacchio2011fast, scacchi2011multilevel, zampini2014dual, huynh2021scalable}. Other works have considered monolithic time discretizations \cite{murillo2004fully, huynh2021newton}, where the resulting algebraic system also includes the ODEs. Implicit-explicit (IMEX) schemes have been also considered, where the nonlinearity is computed using the previous time step \cite{munteanu2009sisc, scacchi2011multilevel, zampini2014dual}.

In all of the aforementioned formulations, a large linear system needs to be solved at each time step. The development of adequate preconditioners for such a system has been successfully addressed with domain decomposition methods \cite{scacchi2011multilevel, huynh2022parallel}. Despite this, all approaches considered thus far inevitably require to perform at least one matrix assembly and one linear solve at each time step, which are computationally expensive. If such a procedure is done only once per time step, as in IMEX schemes, then restrictive time steps are required to guarantee that undesirable numerical errors do not occur, such as oscillation and dispersion. On the other hand, if the solver is embedded inside a Newton nonlinear solver, the timestep solver converges quadratically and more robustly, but at a great computational cost. These two aspects are in general balanced in a highly non-trivial and problem-dependent manner. 

Despite Newton methods being a gold-standard in the numerical approximation of nonlinear PDEs, other methods provide alternative ways of solving the same problem with a hopefully reduced overall complexity. Some of these include modifications of the linearized inverse operator in Newton's method, resulting in quasi-Newton \cite{wright1999numerical} and inexact-Newton \cite{eisenstat1994globally, eisenstat1996choosing} methods, which in general provide superlinear convergence. Other works consider nonlinear variants of well-established linear solvers such as nonlinear Conjugate Gradient (NCG) or nonlinear Generalized Minimal Residual (NGMRES) \cite{washio1997krylov, wright1999numerical, brune2015composing} methods, presenting linear convergence but requiring only the assembly of the residual, without the solution of a linear system. In the context of computational cardiology, such alternatives have been thoroughly studied and compared in cardiac mechanics \cite{barnafi2022parallel} and, preliminarly, in electrophysiology \cite{barnafi2022alternative}. In the latter, the performance of the solvers was initially optimized by varying relevant parameters of each method, and then such optimal configurations were tested in large scale simulations. In all cases, the advantage over a Newton iteration is evident in all of the considered methods. A nonlinear GMRES method for the Bidomain model has also been studied in Ref. \cite{bourgault2003simulation}, and quasi-Newton solvers for nonlinear elasticity have been studied in Refs. \cite{linge2005solving, liu2017quasi}.

In this work, we turn our attention on the analytic study of the convergence properties of these nonlinear solvers when employed for the solution of nonlinear problems arising in staggered solution strategies for cardiac electrophysiology. More specifically, we focus on nonlinear CG method and quasi-Newton methods through the construction of an adequate potential that arises from the ODE/PDE decoupling. We also explore numerically the performance of other methods, such as nonlinear GMRES and inexact-Newton methods to devise nonlinear solvers that are tailored for different physically relevant scenarios. 

Our analysis is based on the construction of a suitable potential, whose first order conditions are given by the time-discretized Bidomain equations. This potential is a particular case of the electrochemical potential considered in Refs. \cite{hurtado2014gradient, cornejo2015analisis}, where the authors show that, whenever such potential exists, it is possible to establish the gradient flow structure of the Bidomain equations and thus guarantee the existence of a unique solution. In particular, the convexity of the potential allows for the computation of an optimal time step, depending only on the problem parameters, that guarantees the convergence of a nonlinear solver.

The work is structured as follows: we first provide a mathematical description of the model in Sec. \ref{sec: Bidomain}, along with the discretization choices and the definition of the functional that will be analyzed next.
A brief overview of the two classes of nonlinear solvers on which our attention is focused, is given in Section \ref{sec: nonlinear solvers}. The main results of the paper are presented in Sec. \ref{sec: convergence analysis}, where we validate all the assumptions and hypotheses needed in order to guarantee and prove the convergence of the methods.
Extensive parallel numerical tests, shown in Sec. \ref{sec: numerical tests}, validate the theoretical results; Section \ref{sec: conclusions} concludes the work by providing closing remarks.

\paragraph{\bf Notation} Throughout this work we will employ basic functional notions, which we report here for a more comprehensive readability. We indicate with $C^s(\omega)$ the space of continuous functions $f:\omega \to \R$ with continuous first $s$ derivatives, and with $H^s(\omega)$ the Hilbert space of functions $f:\omega \to \R$ with norm $\|\cdot \|_{s,\omega}^2\coloneqq (\cdot, \cdot)_{s,\omega}$ and relative seminorm $| \cdot |_{s, \omega}^2$ where we denote the case $s=0$ with $L^2(\Omega)$. For a functional space $V$, $(V)'$ denotes the functional dual space with its norm given by $\sup_{\|u\|_V=1}\langle Vu, u\rangle$, and finally we denote scalar, vector and tensor quantities as $a$, $\vec a$, and $\ten a$ respectively.

\section{The Bidomain Equations for Cardiac Electrophysiology} \label{sec: Bidomain}
We consider the cardiac Bidomain model \cite{franzone2014mathematical}, a system of degenerate parabolic partial differential equations modeling the propagation of the electric signal in the cardiac tissue, known as myocardium. Cardiac tissue can be described electrically as the composition of two ohmic conducting media, the intra- and extracellular domains, separated by the active cellular membrane which acts as insulator between the two domains. This property is fundamental, as otherwise there would be no potential difference across the membrane. In the Bidomain model these anisotropic continuous media are assumed to coexist at every point of the tissue and to be connected by a distributed continuous cellular membrane \cite{franzone2014mathematical}. In this way, we define the electric potential in each point of the two domains as a quantity averaged over a small volume: consequently, we assume that every point of the cardiac tissue belongs to both intracellular and extracellular spaces, thus being assigned both an intra- and an extracellular potential. We will denote by $\Omega$ the cardiac tissue volume represented by the superposition of these two spaces, and in general use the subscripts $i$ and $e$ for intracellular and extracellular quantities respectively. 

The cardiac muscle fibers are modeled as laminar sheets running radially from the outer (epicardium) to the inner surface (endocardium) of the heart, direction in which they present a counterclockwise rotation. Therefore, it is possible to mathematically define the electric conductivity tensors as follows: at each point $\xx$ of the cardiac domain we define an orthonormal triplet of vectors $\aal(\xx)$ parallel to the local fiber direction, $\aat(\xx)$ tangent and orthogonal to the laminar sheets, and $\aan(\xx)$ transversal to the fiber axis \cite{legrice1995laminar}. We define the conductivity tensors $\ten D_i$ and $\ten D_e$ of the two media as
$$
\ten D_{i,e} (\xx) = \sum_{\ast = \left\{ l,t,n \right\} } \sigma_\ast^{i,e} \vec a_\ast (\xx) \vec a_\ast^T(\xx) ,
$$
where $\sigma_{l, t, n}^{i,e}$ are the conductivity coefficients in the intra- and extracellular domains along the corresponding directions. 

In this work we consider the parabolic-parabolic formulation of the Bidomain, which reads: given $I_\text{app}^i, I_\text{app}^e : \Omega \times \mathbb (0,T) \rightarrow \mathbb R$, find the intracellular and extracellular potentials $u_{\star}$: $\Omega \times (0, T) \rightarrow \mathbb{R}$, $\star \in \{i,e\}$, the ionic concentration variables $\vec{c}$: $\Omega \times (0, T) \rightarrow \mathbb{R}^{N^C}$ and the gating variables $\vec{w}$: $\Omega \times (0, T) \rightarrow \mathbb{R}^{N^W}$ (which model the opening and closing process of ionic channels), for $N^C, N^W \in \mathbb{N}$, such that
\begin{equation}\label{eq:Bidomain}
	\begin{dcases}
		\chi C_m \dfrac{\partial v}{\partial t} - \dive (\ten D_i\grad u_i) + \chi I_\text{ion}(v, \vec w, \vec c) = I_\text{app}^i, \\
		-\chi C_m \dfrac{\partial v}{\partial t} - \dive (\ten D_e\grad u_e) - \chi I_\text{ion}(v, \vec w, \vec c) = I_\text{app}^e, \\
		\dfrac{\partial \vec c}{\partial t} - C(v, \vec w, \vec c) = 0,  \qquad
		\dfrac{\partial \vec w}{\partial t} - R(v, \vec w) = 0,  
	\end{dcases}
\end{equation}
with the homogeneous Neumann boundary conditions (i.e. assuming the heart electrically insulated), 
$$
(\ten D_i \grad u_i) \cdot \mathbf n = 0, \qquad
(\ten D_e \grad u_e) \cdot \mathbf n = 0, \qquad
\text{on } \Omega \times (0,T),
$$
where $v = u_i - u_e$ is the transmembrane potential, $C_m$ is the membrane capacitance for unit area of the membrane surface and $\chi$ is the membrane surface to volume ratio. Here $I_\text{app}^{i,e}$ represent the intra- and extracellular applied currents (needed to trigger a propagating front) and initial values 
$$
v(\xx,0) = v_0(\xx), \quad \vec{w}(\xx, 0) = \vec{w}_0 (\xx), \quad \vec{c} (\xx, 0) = \vec{c}_0 (\xx) 				 \qquad \qquad \text{in } \Omega  .
$$
The nonlinear reaction term $I_\text{ion}$ and the functions $C(\cdot, \cdot,\cdot)$ and $R(\cdot, \cdot)$ in the ODEs system for the ionic and gating concentration variables are given by the chosen ionic membrane model. \\
Results on existence, uniqueness and regularity of the solution of system (\ref{eq:Bidomain}) have been extensively studied, see for example Refs. \cite{franzone2014mathematical, veneroni2009reaction, cornejo2015analisis}. 
We recall that, to guarantee the existence of the solution, the following condition must hold
$$
\int_\Omega ( I_\text{app}^i + I_\text{app}^e) \, dx = 0.
$$
Moreover, since the potentials $u_i$ and $u_e$ are unique only up to an arbitrary time dependent constant, in order to fix such constant we impose the condition 
$$
\int_\Omega u_e \, dx = 0.
$$
We highlight that, to our knowledge, the analysis performed in Ref. \cite{cornejo2015analisis} is among the first to study the Bidomain equations in the context of convex analysis. From that point of view, the existence of solutions depends on the existence of a gradient flow formulation of problem \eqref{eq:Bidomain}, which at the same time depends on the existence of an electrochemical potential $F$ such that
$$ \frac{\partial F}{\partial v} = I_\text{ion}, \quad \frac{\partial F}{\partial \vec w} = R, \quad \frac{\partial F}{\partial \vec c} = C, $$
which is highly non trivial, and in some cases possibly not even true. In Ref. \cite{cornejo2015analisis}, the construction of $F$ is performed in the simple case of the FitzHugh-Nagumo ionic model. Even if the existence of an electrochemical potential for any triplet $(I_\text{ion}, \, R, \, C)$ is out of our scope, in what follows it is important to notice that a potential for $I_\text{ion}$ can be constructed by means of integration, whenever $\vec w$ and $\vec c$ are fixed, and this automatically grants a variational structure to the PDEs in \eqref{eq:Bidomain}.

\subsection{Time and space discretizations}\label{sec: discretization}
We consider standard first order finite elements for the discretization of (\ref{eq:Bidomain}). Since the subsequent analysis is independent from the space discretization choice, we will not give any details of it, but the interested reader can refer to Ref. \cite{franzone2014mathematical}. 

For our analysis, we consider an ODE/PDE decoupling solution strategy for problem (\ref{eq:Bidomain}), in the same fashion as Refs. \cite{huynh2022parallel, munteanu2009decoupled, zampini2014dual}, where the gating and ionic concentration variables $\vec w, \vec c$ are solved for a given previous transmembrane potential $v^n\approx v(t^n)$, considered in a time discrete scenario. Thus, at each time step, the ODEs system representing the ionic model is solved first; then, the nonlinear algebraic Bidomain system is solved and updated.
In a very schematic way, this decoupling strategy can be summarized as follows: for each time step $n$, 

\begin{enumerate}
	\item Given the intra- and extracellular potentials at the previous time step, define $v^{n-1} \coloneqq u_i^{n-1} - u_e^{n-1}$ and compute the gating and ionic concentrations variables $\vec w^n, \vec c^n$ such that
	$$
		\dfrac{\vec c^n - \vec c^{n-1}}{\tau} + C(v^{n-1}, \vec w^n, \vec c^n) = 0, 
  \qquad
		\dfrac{\vec w^n - \vec w^{n-1}}{\tau} + R(v^{n-1}, \vec w^n) = 0.
	$$
	\item Solve and update the Bidomain nonlinear system. \\Given $u_{i,e}^{n-1}$ at the previous time step and given $\vec w^{n}$ and $\vec c^{n}$ (from step 1), compute $\vec u^{n} = (u_i^{n}, u_e^{n})$ by solving the nonlinear system
	\begin{equation}\label{eq:Bidomain-discrete}
		\begin{cases}
			\chi C_m \dfrac{v^n - v^{n-1}}{\tau} - \dive \ten D_i\grad u_i^n + \chi I_\text{ion}(v^n, \vec w^n, \vec c^n) &= I_\text{app}^i,\\
			-\chi C_m \dfrac{v^n - v^{n-1}}{\tau} - \dive \ten D_e\grad u_e^n - \chi I_\text{ion}(v^n, \vec w^n, \vec c^n) &= I_\text{app}^e. 
		\end{cases}
	\end{equation}
\end{enumerate}
We denote with $\tau = t^n - t^{n-1}$ the timestep. This strategy is usually adopted in contrast to a monolithic approach, where the two systems are solved together, and where the computational workload is higher due to the presence of the ionic model in the nonlinear algebraic system. The ODE/PDE decoupling strategy presents the advantage of allowing for different time scales for the different systems of equations, and it has been extensively studied together with several scalable parallel preconditioners, e.g. Refs. \cite{huynh2021newton, munteanu2009sisc, murillo2004fully}. 
We remark that this decoupling approach does not impair the order of accuracy of the method, which in this case remains of order 1 in time.
\rev{We also note that, differently from the more popular IMEX methods, our strategy consists of treating the nonlinear reaction term
implicitly, i.e. putting $v^n$ instead of $v^{n-1}$. This yields at each time step the solution of a nonlinear algebraic system.
Due to the nonlinearity of the reaction term, even in our strategy the time step size is subject to a stability constraint that
depends on the derivative of the ionic term, see \cite{marsh2012secrets}.
However, this constraint should be milder than using an IMEX method, thus allowing slightly larger time steps.}
In what follows, we focus on the robust and efficient solution of the nonlinear system \eqref{eq:Bidomain-discrete}.

\begin{remark}
    Many other numerical integration schemes are available in the literature that are suitable for system \eqref{eq:Bidomain}. Our analysis is independent of the discretization considered for the ODE system, and depends only on (i) having a splitting approach that decouples the ODEs and PDEs and (ii) having an implicit time discretization of the PDE. On this line, one simple application of this work would be using a higher order integration for the ODEs and a Crank-Nicholson for the PDE.
\end{remark}

\subsection{Finding a suitable time-discrete Bidomain potential}\label{sec: potential}
In virtue of the decoupling strategy described in Section \ref{sec: discretization}, we focus only on equation \eqref{eq:Bidomain-discrete}. We consider the partial primitive of $I_\text{ion}$ as follows,
\begin{equation} \label{eq:theta}
	\Theta(v,\vec w) = \int_{v_0}^v \chi I_\text{ion}(\xi, \vec w)\,d\xi,
\end{equation}
where, without loss of generality, we will only write $\vec w$ instead of $\vec w,\vec c$. This yields
$$
\partial_v \Theta(v,\vec w) = \chi I_\text{ion}(v, \vec w) \quad \forall \vec  w\in\mathbb R^{N^W},
$$
where we denote the partial derivative with respect to $v$ as $\partial_v \coloneqq \frac{\partial}{\partial v}$. In particular it will hold that, for $s>0$, if $I_\text{ion}(\cdot, \vec w)\in H^s(\R)$ for all $\vec w$, then $\Theta(\cdot, \vec w)\in H^{s+1}(\R)$ for all $\vec w$. This also holds for strong continuity, i.e. in $C^s(\R)$.

By defining the spaces $V=H^1(\Omega)$ and $\widetilde{V}=\{ \mu\in V: \int_\Omega \mu\,dx = 0\}$, it is possible to define the Bidomain potential $\Psi:V\times \widetilde{V}\mapsto \R$ as
\begin{multline}\label{eq:potential}
	\Psi(u_i, u_e) = \frac 1 2 \int_\Omega \frac{\chi C_m}{\tau} (v-v^{n-1})^2 \,dx + \frac 1 2 \int_\Omega \left(\ten D_i\grad u_i\right)\cdot\grad u_i\,dx + \frac 1 2 \int_\Omega \left(\ten D_e\grad u_e\right) \cdot \grad u_e \,dx  \\
	+ \int_\Omega \Theta(u_i - u_e, \vec w^n) \,dx - \int_\Omega I_\text{app}^i u_i\,dx - \int_\Omega I_\text{app}^e u_e \,dx,  
\end{multline}
whose stationary points are given by the weak formulation of \eqref{eq:Bidomain-discrete}. Indeed, this can be easily verified by computing the partial Gateaux derivatives of $\Psi$:
\begin{multline*}
	\partial_{u_i}\Psi(u_i, u_e)[\varphi_i] \coloneqq  \frac{d}{d\epsilon}\bigg|_{\epsilon=0}\Psi(u_i + \epsilon \varphi_i, u_e) \\
	= \int_\Omega \frac{\chi C_m}{\tau} (v-v^{n-1}) \varphi_i\,dx + \int_\Omega \left(\ten D_i\grad u_i\right)\cdot\grad \varphi_i\,dx + \int_\Omega \chi I_\text{ion}(v, \vec w^n)\,\varphi_i\,dx - \int_\Omega I_\text{app}^i \varphi_i\,dx,
\end{multline*}
\begin{multline*}
	\partial_{u_e}\Psi(u_i, u_e)[\varphi_e] \coloneqq \frac{d}{d\epsilon}\bigg|_{\epsilon=0}\Psi(u_i, u_e + \epsilon \varphi_e) \\
	= - \int_\Omega \frac{\chi C_m}{\tau} (v-v^{n-1}) \varphi_e\,dx + \int_\Omega \left(\ten D_e\grad u_e\right)\cdot\grad \varphi_e\,dx - \int_\Omega \chi I_\text{ion}(v, \vec w^n)\,\varphi_e\,dx - \int_\Omega I_\text{app}^e \varphi_e\,dx.
\end{multline*}

\noindent
Note that the changed signs in the time derivative and in $I_\text{ion}(v, \vec w^n)$ arise naturally from the derivation of $v$ due to the chain rule. Therefore, we can formulate at each time instant $t^n$ problem \eqref{eq:Bidomain-discrete} as the solution of the following minimization problem: 
\begin{equation} \label{eq: minimization problem}
	(u_i^n, u_e^n) = \argmin_{(u_i, u_e)\in V\times\widetilde{V}} \, \Psi(u_i, u_e).
\end{equation}
For simplicity, from now on, we drop the index $n$, if ambiguity does not occur. \rev{We also highlight that we refer to this approach as implicit, due to the treatment of the ionic current term. Still, it is a semi-implicit scheme regarding the coupled ODE-PDE system.}

\section{Nonlinear Bidomain Solvers} \label{sec: nonlinear solvers}
In this section, we provide a review of the nonlinear Bidomain solvers under consideration, together with their convergence theory. All of the following methods could be extended with a step for computing the step length.
This extension is beyond the scope of this paper and use only residual based solvers implemented in the SNES package of PETSc \cite{balay2019petsc} with a full step length of $\alpha_k=1$. Details on step length computation algorithms can be found in \cite{wright1999numerical}.

\subsection{Quasi-Newton (QN) methods } \label{sec: quasi-newton}
Quasi-Newton methods consider a simplified Newton step, where the Jacobian is never computed but is approximated at each iteration. QN methods require only the gradient of the function to minimize to be provided at each iteration. We report here a brief and not exhaustive overview of a subclass of QN algorithms, the BFGS method (named after Broyden, Fletcher, Goldfarb and Shanno); for more details we refer to Ref. \cite{wright1999numerical}.\\

\noindent
Consider the minimization problem $\min_x f(x)$ and its quadratic model at the current iterate $x_k$, $$\min_p (f_k + \nabla f^T_k p + \frac{1}{2} p^T B_k p)$$ with $B_k$ a symmetric positive definite matrix. 
The minimizer $p_k = -B^{-1}_k \nabla f_k$ is used as search direction for the update of the iteration $x_{k+1} = x_k + \alpha_k p_k$, where $\alpha_k$ is a suitable step length; see Algorithm \ref{alg: bfgs} below.

This iteration resembles the line search Newton, though the main difference is that $B_k$ approximates the true Hessian matrix, along with an appropriate initialization $B_0$. Then, at each subsequent iteration, this approximation is enriched with the previous iterations by means of the following BFGS updating formula
\begin{equation} \label{eq: approx hessian}
	B_{k+1}^{-1} = (I - \rho_k s_k y_k^T) B_k^{-1} (I - \rho_k y_k s_k^T) + \rho_k s_k s_k^T,
\end{equation}
where $s_k = x_{k+1} - x_k$, $y_k = \nabla f(x_{k+1}) - \nabla f(x_k)$ and $\rho_k = 1/\langle s_k, y_k \rangle$.

\begin{algorithm}
	\caption{Quasi-Newton method, BFGS algorithm}\label{alg: bfgs}
	\begin{algorithmic}[1]
		\State given initial guess $x_0$ and convergence tolerance $\epsilon > 0$
		\State compute the approximation of the inverse of the Jacobian $B_0^{-1}$
		\State $k \gets 0$
		\While{$\| \nabla f_k \| > \epsilon$}
		\State compute search direction $p_k = - B_k^{-1} \nabla f_k$
		\State compute step length $\alpha_k$ (in our case, $\alpha_k=1$) and set $x_{k+1} = x_k + \alpha_k p_k$
		\State compute $s_k, y_k$ and $\rho_k$
		\State compute the action of $B_{k+1}^{-1}$ by means of Eq. (\ref{eq: approx hessian})
		\State $k \gets k+1$
		\EndWhile
	\end{algorithmic}
\end{algorithm}

Each iteration of Algorithm \ref{alg: bfgs} can be performed with a cost of $O(n)$ operations, since it does not require the solution of linear systems or matrix-matrix operations. Its rate of convergence can be proved to be superlinear (see Ref. \cite{wright1999numerical}) which, of course, makes it converge less rapidly than Newton's method (which converges quadratically), but with a greatly reduced computational cost per iteration, since there is no need for the second derivative, and the action of $B_k$ is given by a recursive formula that hinges only on the inversion of $B_0$.

The convergence of the method is guaranteed under the conditions described in Ref. \cite[Theorem 6.6]{wright1999numerical}, which we have adapted to our setting.

\begin{theorem}[BFGS convergence]\label{thm:bfgs}
	Consider the level set of the potential $\Psi$ defined in (\ref{eq:potential})
	$ \mathcal L(\hat{\vec u}) \coloneqq \{\vec u=(u_i, u_e) \in V \times \widetilde V: \Psi(\vec u)\leq \Psi(\hat{\vec u})\}$
	together with an initial guess $\vec u_0 \in V\times \widetilde{V}$ and the following properties:
	\begin{itemize}
		\item [\bf A1] the objective function $\Psi$ is twice continuously differentiable, i.e. the partial derivatives $\partial_{u_ju_k}\Psi: V\times \widetilde{V}\to (V\times \widetilde{V})'$, for $\{j,k\}\in\{i,e\}$ are continuous in the corresponding norms: 
		$$ \|\vec u - \hat{\vec u} \|_{V\times \widetilde{V}}\to 0 \quad\Rightarrow\quad \| \partial_{u_ju_k}\Psi(\vec u) - \partial_{u_ju_k}\Psi(\hat{\vec u})\|_{(V\times \widetilde{V})'}\to 0$$
		for all $\vec u, \hat{\vec u} \in V\times \widetilde{V}$.
		\item [\bf A2] the level set $\mathcal L(\vec u_0)$ is convex and there exist positive constants $m$ and $M$ such that 
		$$
		m \| \vec z \|^2_{V} \leq d^2\Psi(\vec u)[\vec z] \leq M \| \vec z \|^2_{V},
		$$
		for all $\vec z \in V$ and $\vec u \in \mathcal L(\vec u_0)$, where $d^2\Psi$ is the second variation (Hessian) of $\Psi$;
		\item [\bf A3] the Hessian $d^2 \Psi$ is Lipschitz continuous at the minimum $\vec u^\star$, that is
		$$
		\| d^2 \Psi (\vec u) - d^2 \Psi (\vec u^\star) \|_{(V\times \widetilde{V})'} \leq L \| \vec u - \vec u^\star \|_{V\times \widetilde{V}} ,
		$$
		for all $\vec u$, where $L$ is a positive constant.
	\end{itemize}
	If all the above properties hold, the BFGS algorithm converges superlinearly. 
\end{theorem}

\subsection{Nonlinear Conjugate Gradient} \label{sec: nonlinear CG}
It is well established that the Conjugate Gradient (CG) method for the iterative solution of a generic linear system $Ax-b$ can be reformulated as a minimization algorithm for the convex quadratic function $$\phi(x) = \frac{1}{2} x^TAx - b^Tx.$$ Fletcher and Reeves \cite{fletcher1964function} shows how to modify the linear CG in order to adapt this approach to minimize general nonlinear functions, resulting in Algorithm \ref{alg: ncg}.

\begin{algorithm}[ht!]
	\caption{Nonlinear CG method, Fletcher -- Reeves algorithm (NCG-FR)}\label{alg: ncg}
	\begin{algorithmic}[1]
		\State given initial guess $x_0$;
		\State evaluate $f_0 = f(x_0)$, $\nabla f_0 = \nabla f(x_0)$
		\State $p_0 \gets - \nabla f_0$ and $k \gets 0$
		\While{$\nabla f_k \neq 0$}
		\State compute step length $\alpha_k$ (in our case, $\alpha_k=1$) and set $x_{k+1} = x_k + \alpha_k p_k$
		\State evaluate $\nabla f_{k+1}$:
		\begin{align}
			\beta_{k+1} &\gets \dfrac{\nabla f^T_{k+1} \nabla f_{k+1}}{\nabla f^T_k \nabla f_k} \\
			p_{k+1} &\gets - \nabla f_{k+1} + \beta_{k+1} p_k \label{eq: ncg search direction}\\    
			k &\gets k+1
		\end{align}
		\EndWhile
	\end{algorithmic}
\end{algorithm}

Other variants of the algorithm differ mainly on the computation of $\beta_{k+1}$. Additionally, a common modification adopted in the numerical implementation of Algorithm \ref{alg: ncg} is the introduction of a restart. This means that every $m$ iterations $\beta_k = 0$, discarding old information that may not be useful anymore and refreshing the algorithm.
As pointed out in Ref. \cite{wright1999numerical}, in practical context where $m$ is usually large, restart may never occur, since the approximate solution may be obtained with fewer iterations. In this sense, numerical implementation of nonlinear CG may present different strategies for restarting.

The convergence of the method is guaranteed under the conditions shown in Ref. \cite[Theorem 5.6]{wright1999numerical}. 

\begin{theorem}[NCG convergence]
	Consider the level set
	$ \mathcal L(\vec u) \coloneqq \{\vec w: \Psi(\vec w)\leq \Psi(\vec u)\}$
	together with an initial guess $\vec u_0$ and the following properties:
	\begin{itemize}
		\item[\bf B1.] The level set $\mathcal L(\vec u_0)$ is bounded.
		\item[\bf B2.] The objective function $\Psi$ is Lipschitz differentiable.
	\end{itemize}
	If both conditions hold, then the NCG algorithm converges linearly.
\end{theorem}

\subsubsection{Mesh-independence property}

One of our main concerns when developing a solver is its optimality, understood as being robust with respect to number of degrees of freedom. In the literature, this concept is known as the mesh-independence property: we briefly report here two approaches, which are reformulations of the same methods but on an infinite-dimensional setting. We highlight that so far these theories are well-understood only for Newton and Quasi-Newton methods.

\paragraph{\textbf{Quasi-Newton methods}} In addition to the hypotheses shown in Theorem \ref{thm:bfgs}, the initial approximation of the Hessian operator must differ from it up to a compact operator \cite{sachs1986broyden, griewank1987}. This hypothesis has not been further studied, and the existence of more verifiable theoretical framework for its use remains an open problem. This in particular means that the mesh-independence of Quasi-Newton methods for the Bidomain equations might not hold if the initial approximation of the operator is not sufficiently good.

\paragraph{\textbf{Newton methods}} In \cite{weiser2005asymptotic} the authors verify that, using an affine-invariant property of the Jacobian, it is possible to establish the asymptotic mesh-independence of Newton methods for a large class of equations. In particular, they investigate this property for abstract elliptic semilinear equations, meaning that, since the Bidomain equations belong to this category of equations, these are expected to display the mesh-independence property whenever solved with a Newton method.  \\

In addition to these theories, we note that such a theory for first order methods has not been established yet in case of non-quadratic nonlinear problems, and indeed the numerical tests may present unpredictable behaviors (see Section \ref{sec: numerical tests}).

\section{Convergence Analysis for the Bidomain Equations} \label{sec: convergence analysis}
In this section we first provide several assumptions required for proving the convergence of the BFGS and Fletcher--Reeves algorithms for the Bidomain equations \eqref{eq:Bidomain-discrete}. After this, we verify that the proposed assumptions yield the corresponding properties of each method.

The assumptions under consideration are the following: 
\begin{itemize}
	\item {\bf (H1)} The function $\Theta$ defined in (\ref{eq:theta}) is bounded from below and grows with a prescribed rate $\alpha\geq 1$:
	$$ \Theta(x, \vec w) \geq c_1|x|^\alpha + \theta_0,$$
	where $\theta_0\in \R$ and $c_1>0$ may depend on $\vec w$.
	\item {\bf (H2)} The function $\Theta(\cdot, \vec w)$ is Lipschitz differentiable, i.e. its derivative is Lipschitz continuous. 
	\item {\bf (H2*)} The function $\Theta(\cdot, \vec w)$ is twice Lipschitz differentiable, i.e. its first and second order derivatives are Lipschitz continuous. 
	\item {\bf (H3)} The function $\partial_v I_\text{ion}(v,\vec w)$ is uniformly bounded from by two constants $\underline I \leq \overline I \in \R$: for all $v,\vec w$, it holds\footnote{
        We underline that this assumption cannot be theoretically proved, since for the Bidomain equations it does not exists a maximum principle; however, since usually $v$ and $\vec w$ belong to a fixed range of values, this assumption holds from a numerical viewpoint.
    }
	$$\underline I \leq \partial_v I_\text{ion}(v,\vec w) \leq \overline I. $$
	\item {\bf (H4)} The conductivity tensors are symmetric and positive definite: there exist positive constants $D_i^0,D_e^0 \in \R$ such that
	$$ \int_\Omega \grad v^T\ten D_\ast\grad v\,dx\geq D_\ast^0\normL{\grad v}^2, \quad\forall v\in H^1(\Omega), \quad \ast\in \{i,e\}.$$
    From now on, in particular, we will refer to $D^0 \coloneqq \min \{ D^0_i, D^0_e \}$. 
\end{itemize}

We consider separately \textbf{(H2)} and \textbf{(H2*)} due to the different regularity conditions of both methods under consideration. We note that the assumptions in particular guarantee that $\Psi$ has a minimum. This can be seen from the coercivity given by \textbf{(H1)} and \textbf{(H4)}, plus the continuity from \textbf{(H2)}. The conclusion is drawn using the direct method of the calculus of variations \cite{dacorogna2007direct, gelfand2000calculus}. We note also that {\bf (H1)} is a requirement in the convergence of many optimization methods and, at a continuous level, it is the one that guarantees the boundedness of the infimum, which is the point of departure for more sophisticated techniques in convex analysis.

Assumption  \textbf{(H3)} is the most restrictive, and it depends on the choice of $I_\text{ion}(\cdot, \cdot)$: in general, given a ionic model, it can be checked numerically if this assumption holds for the values under consideration, and indeed as our numerical results show, this Assumption is not limiting in practice. Finally, we observe that \textbf{(H4)} may not hold in pathological scenarios, due to areas with no conductivity, yielding positive semi-definite conductivity tensors.

\subsection{Convergence analysis of Quasi-Newton methods (BFGS algorithm)} \label{sec: convergence BFGS}

We now prove, for the Bidomain setting, the three conditions needed for the convergence of quasi-Newton methods. 

\paragraph{\bf Property A1 (twice differentiability of the objective function)} Proving that the objective function $\Psi$ is twice continuously differentiable means to prove that
$$
\partial^2_{u_i u_i} \Psi (u_i, u_e),
\qquad
\partial^2_{u_i u_e} \Psi (u_i, u_e),
\qquad
\partial^2_{u_e u_e} \Psi (u_i, u_e),
$$
exist and are continuous. Since $v = u_i - u_e$ and thanks to the chain rule, we have
$$
\partial_{u_i} \, g(v, \vec w) = \frac{dv}{du_i} \, \frac{d}{dv} \, g(v, \vec w) = \partial_v \, g(v, \vec w) = -\frac{dv}{du_e} \, \frac{d}{dv} \, g(v, \vec w) = -\partial_{u_e} \, g(v, \vec w).
$$
By computing the partial Gateaux derivatives of $\Psi$ and using the dominated convergence theorem\footnote{
	Indeed, in virtue of assumption \textbf{(H2)}, we have that 
	$$ \frac{I_\text{ion}(v-\epsilon \phi_e, \vec w) - I_\text{ion}(v,\vec w)}{\epsilon}\leq L_{I_\text{ion}}|\phi_e|,$$
	where $L_{I_\text{ion}}$ is the Lipschitz continuity constant of $I_\text{ion}$, and $\phi_e(x)$ is measurable. }  in virtue of \textbf{(H2*)}, we obtain
\begin{align*}
	&\partial^2_{u_i u_i}\Psi(u_i, u_e) [\varphi_i, \phi_i] = \frac{d}{d\epsilon}\bigg|_{\epsilon=0} \partial_{u_i}\Psi(u_i + \epsilon \phi_i, u_e)[\varphi_i]  \\
	&= \int_\Omega \frac{\chi C_m}{\tau} \varphi_i \phi_i \,dx + \int_\Omega \left(\ten D_i\grad \varphi_i\right)\cdot\grad \phi_i\,dx + \lim_{\epsilon \rightarrow 0} \int_\Omega \chi \, \frac{ I_\text{ion}(v + \epsilon \phi_i, \vec w) - I_\text{ion} (v, \vec w) } {\epsilon} \,\varphi_i\,dx \\
	&= \int_\Omega \frac{\chi C_m}{\tau} \varphi_i \phi_i \,dx + \int_\Omega \left(\ten D_i\grad \varphi_i\right)\cdot\grad \phi_i\,dx + \int_\Omega \chi \partial_v I_\text{ion} (v,\vec w) \,\varphi_i \phi_i\,dx ,
\end{align*}
\begin{align*}
	&\partial^2_{u_e u_e}\Psi(u_i, u_e)[\varphi_e, \phi_e] = \frac{d}{d\epsilon}\bigg|_{\epsilon=0} \partial_{u_e}\Psi(u_i, u_e + \epsilon \phi_e)[\varphi_e]  \\
	&= \int_\Omega \frac{\chi C_m}{\tau} \varphi_e \phi_e \,dx + \int_\Omega \left(\ten D_e \grad \varphi_e\right)\cdot\grad \phi_e\,dx - \lim_{\epsilon \rightarrow 0} \int_\Omega \chi \, \frac{ I_\text{ion}(v - \epsilon \phi_e, \vec w) - I_\text{ion} (v, \vec w) } {\epsilon} \,\varphi_e\,dx \\
	&= \int_\Omega \frac{\chi C_m}{\tau} \varphi_e \phi_e \,dx + \int_\Omega \left(\ten D_e\grad \varphi_e\right)\cdot\grad \phi_e\,dx + \int_\Omega \chi \partial_v I_\text{ion} (v, \vec w) \,\varphi_e \phi_e \,dx  ,
\end{align*}
\begin{align*}
	\partial^2_{u_i u_e}\Psi(u_i, u_e)&[\varphi_i, \phi_e] = \frac{d}{d\epsilon}\bigg|_{\epsilon=0} \partial_{u_i}\Psi(u_i , u_e +  \epsilon \phi_e)[\varphi_i]  \\
	&= \int_\Omega \frac{\chi C_m}{\tau} \varphi_i \phi_e \,dx + \lim_{\epsilon \rightarrow 0} \int_\Omega \chi \, \frac{ I_\text{ion}(v - \epsilon \phi_e, \vec w) - I_\text{ion} (v, \vec w) } {\epsilon} \,\varphi_i\,dx \\
	&= \int_\Omega \frac{\chi C_m}{\tau} \varphi_i \phi_i \,dx - \int_\Omega \chi \partial_v I_\text{ion} (v, \vec w) \,\varphi_i \phi_e\,dx .
\end{align*}
We do not compute the other crossed derivative as they match due to symmetry of the potential $\Psi$.
These second derivatives are continuous if and only if, for a given $\vec w$, the functional 
\begin{equation} \label{eq:integral-Dv-Iion}
	\mathcal I(\vec u) = \int_{ \Omega} \chi \partial_v I_\text{ion} (v, \vec w) \varphi \phi dx, 
	\qquad
	\forall \varphi, \phi \in \Hone
\end{equation}
is continuous in the sense of \textbf{(A1)}. Since it is possible to bound the right-hand side of (\ref{eq:integral-Dv-Iion}) with 
$$
C_\Omega \, \max_{v} | \partial_v I_\text{ion} (v,\vec w) | \, \| \varphi \|_{H^1(\Omega)} \, \| \phi \|_{H^1(\Omega)},
$$
with $C_\Omega$ positive constant depending on the domain $\Omega$, we only have to work on $\partial_v I_\text{ion} (v,\vec w)$. 
We note that, by fixing $\vec w$, we have $I_\text{ion} (v) \coloneqq I_\text{ion} (v,\vec w)$ such that $\partial_v I_\text{ion} (v,\vec w) = I'_\text{ion}(v)$. Then the conclusion follows from assumption \textbf{(H2*)}.

\paragraph{\bf Property A2 (convexity of the level set)}
In order to prove that the level set $\mathcal L$ is convex, it is sufficient to prove that $\Psi$ is convex. The objective function has already been proven to be twice continuously differentiable. We need now to prove that its second variation $d^2 \Psi (u_i, u_e)[(\phi_i, \phi_e)]$ is positive for $(\phi_i, \phi_e)$ in $V\times \widetilde{V}$. Letting $\vec \phi = (\phi_i, \phi_e)$, we have
\begin{equation}\label{eq:Psi-convexity}
	\begin{aligned}
		&d^2 \, \Psi (u_i, u_e) [(\phi_i, \phi_e)] \\
		&= \int_\Omega \frac{\chi C_m}{\tau} (\phi_i-\phi_e)^2 + \sum_{\ast = i,e} \int_\Omega (\ten D_\ast \nabla \phi_\ast) \cdot \nabla \phi_\ast + \int_\Omega [\phi_i, \phi_e] \, \nabla^2_{(u_i, u_e)} \Theta(v,w) \, [\phi_i, \phi_e]^T \\
		&\geq \frac{\chi C_m}{\tau} \normL{ \phi_i - \phi_e }^2 + \min\{D_i^0, D_e^0\}\normL{\grad {\vec \phi}}^2 + \underline I \normL{\phi_i-\phi_e}^2 \\
		&= \left( \frac{\chi C_m}{\tau} + \underline I  \right) \normL{\phi_i - \phi_e}^2 + \min\{D_i^0, D_e^0\} \normL{\grad {\vec \phi}}^2,
	\end{aligned}
\end{equation}
where we used Assumptions {\bf (H3)} and {\bf (H4)}.
\newline
If $\Theta$ is not convex, then necessarily $\underline I <0$. However, we can impose a restriction on the time step $\tau$ as in Ref. \cite{hurtado2014gradient},
\begin{equation}\label{eq:dt-bound}
	\left( \frac{\chi C_m}{\tau} -|\underline I|\right) \geq 0 \quad\Leftrightarrow\quad \tau \leq \frac{\chi C_m}{|\underline I|}.    
\end{equation}
In this way, we can ensure convexity of the potential $\Psi(u_i, u_e)$ on the bounded domain $\Omega$. 

\begin{remark}
	We highlight that estimate \eqref{eq:dt-bound} is independent of the conductivity tensors. Still, the conductivities contribute to the overall convexity of the potential $\Psi$ as shown in \eqref{eq:Psi-convexity}, which justifies possible degradation of the method in pathological scenarios. 
\end{remark}

\paragraph{\bf Property A3 (Lipschitz continuity of the Hessian matrix)}
This has been proven in Property \textbf{A1} and is a straightforward consequence of assumption \textbf{(H2*)}.

\subsection{Convergence analysis of the NCG-FR method} \label{sec: convergence nonlinear CG}
We proceed by proving the required conditions for NCG-FR, as done for the BFGS method. 

\paragraph{\bf Property B1 (boundedness of the level set)}
We consider a point $\vec u \in \mathcal L$ so that $\Psi(\vec u)\leq \Psi(\vec u_0)$ by definition. We first observe that, by using the inequality $2ab\leq \epsilon a^2 + \epsilon^{-1}b^2$, with $\epsilon$ an arbitrary positive constant, it holds 
$$
\vec I_\text{app}\cdot \vec u = I_\text{app}^i u_i + I_\text{app}^e u_e \leq \frac \epsilon 2 \left((I_\text{app}^i)^2 + (I_\text{app}^e)^2\right) + \frac {\epsilon^{-1}}2\left(u_i^2 + u_e^2\right)
= \frac \epsilon 2 | \vec I_\text{app} |^2 + \frac {\epsilon^{-1} }2 | \vec u |^2,
$$
where $\vec I_\text{app} = (I_\text{app}^i, I_\text{app}^e)$.
Thus, thanks to assumption \textbf{(H1)}, we obtain
$$
\Psi(\vec u_0)\geq \Psi(\vec u) \geq D^0\seminormH{\vec u}^2 + \int_\Omega \theta_0\,dx - \frac \epsilon 2\normL{\vec I_\text{app}}^2 - \frac {\epsilon^{-1}}2\seminormH{\vec u}^2,
$$
and, by rearranging the terms and using assumption \textbf{(H4)}, 
$$ \Psi(\vec u) + \frac \epsilon 2 \normL{\vec I_\text{app}}^2 - \int_\Omega \theta_0\,dx \geq \left(D^0 - \frac {\epsilon^{-1}}2\right)\seminormH{\vec u}^2. $$
Notice that we consider $\epsilon$ such that $\epsilon > \frac 1 {2D^0}$. 
Moreover, we observe that we do not need to control the term involving $\theta_0$, since by definition it is the lower bound of $\Theta$, thus the left-hand side of the above inequality cannot become negative. We verify now that the left-hand side is non-zero:
\begin{align*}
	&\Psi(\vec u) + \frac\epsilon 2\normL{\vec I_\text{app}}^2 - \int_\Omega\theta_0\,dx \\
	&\geq \frac{\chi C_m}{2\tau}\normL{v^0-v^n}^2 + D^0\seminormH{\vec u}^2 + \int_\Omega\Theta(v^0, w)\,dx - \int_\Omega\vec I_\text{app}\cdot \vec u\,dx +\frac \epsilon 2\normL{\vec I_\text{app}}^2\\
	&\geq \frac{\chi C_m}{2\tau}\normL{v^0-v^n}^2 + (D^0 - \frac {\epsilon_2^{-1}}2)\seminormH{\vec u_0}^2 + \int_\Omega c_1|v^0|^\alpha + \frac 1 2(\epsilon - \epsilon_2)\|\vec I_\text{app}\|^2,
\end{align*}
where $v^0$ is the given initial value of the transmembrane potential.
In order for the above inequality to be positive, we have to require that i) $\epsilon_2\geq 1/{2D^0}$ and ii) $\epsilon \geq \epsilon_2$, which holds if e.g. $\epsilon_2 = 1/{2D^0}$. Indeed, i) holds trivially, and by definition we obtain ii) from $\epsilon > 1/{2D^0} = \epsilon_2$. 

We observe that the inequality $\epsilon>\epsilon_2$ is fundamental, since a common initial guess is $v^0=v^n$, from which the left-hand side can be null whenever $\epsilon=\epsilon_2$. 
We conclude that the left hand side is indeed positive and thus $\mathcal L(\vec u_0)$ is bounded.

\paragraph{\bf Property B2 (differentiability of the objective function)} This is shown verbatim as \textbf{A1} from Section \ref{sec: convergence BFGS} by using assumption \textbf{(H2)}.

We summarize the conditions required by each method in Table \ref{table:assumptions}. 
It is interesting to note that only the Fletcher--Reeves algorithm considers assumption \textbf{(H1)}, and instead only BFGS uses assumption \textbf{(H3)}. Indeed, condition \textbf{(H1)} is stronger, as it implies \textbf{(H3)}, which is a reasonable result: Fletcher--Reeves algorithm requires slightly more regular potentials than BFGS in terms of their growth rate. Still, BFGS requires more regularity, as it considers assumption \textbf{(H2*)} instead of \textbf{(H2)}. 
In conclusion, there is no clear answer on which method requires more constraints on the energy function: BFGS requires more regularity, but on the other hand Fletcher--Reeves requires specific growth conditions. Numerical evidences suggest that the numerical approximation is bounded, thus making the Assumptions needed by the Fletcher--Reeves less sharp.

\begin{table}[t]
	\centering
	\caption{List of assumptions (\textbf{H}) required by the quasi-Newton BFGS method (\textbf{A} properties) and by the nonlinear CG, Fletcher--Reeves (NCG-FR) algorithm (\textbf{B} properties).}
	\label{table:assumptions}
	\begin{tabular}{@{}cll@{}}
		\cmidrule[\heavyrulewidth]{2-3}
		& Property & Assumptions                                                       \\
		\midrule
		\multirow{3}{*}{QN-BFGS}    &\textbf{A1} (twice differentiability of the obj. function) & \textbf{(H2*)}       \\
		&\textbf{A2} (convexity of the level set) & \textbf{(H3), (H4)}                                              \\
		&\textbf{A3} (Lipschitz continuity of the Hessian matrix) & \textbf{(H2*)}                                                   \\
		\midrule 
		\multirow{3}{*}{NCG-FR}     &\textbf{B1} (boundedness of the level set) & \textbf{(H1), (H4)}  \\
		&\textbf{B2} (differentiability of the obj. function) & \textbf{(H2)}                                                    \\ 
		\bottomrule
	\end{tabular}
\end{table}

Additionally, we note that the convergence theory of the BFGS holds for the entire Broyden family of methods, and also the convergence of the Fletcher--Reeves yields the same result for many other variants.

\section{Numerical Tests} \label{sec: numerical tests}
In this Section we present the main numerical tests of this work. These consist of: i) studying the robustness of the methods with respect to the problem size, ii) studying the robustness of the methods with respect to discontinuities in the conductivity coefficients modeling the presence of an ischemic region, iii) studying the behavior of the solvers throughout the entire evolution of the propagation of the action potential, iv) studying the parallel scalability of the methods and v) verifying the convergence of the methods. All tests are performed using the optimal parameters obtained in Section \ref{sec:tuning}.

\subsection{Setting}
Unless otherwise stated, we simulate the propagation of the electric signal over a short time interval of $1$ ms, where the initial activation phase of the Bidomain model is computationally most intense. We consider an idealized left ventricle geometry, modeled as a portion of a truncated ellipsoid. We apply an external current $I^i_\text{app}=-I^e_\text{app}=100\,\frac{\text{mA}}{\text{cm}^3}$ for $1$ ms in a small portion of the endocardial surface. 
If not otherwise specified, the Ten Tusscher--Panfilov ionic model \cite{ten2004model} is employed.
The conductivity coefficients values can be found in \cite[Table 6.1]{huynh2022parallel}.
The time interval is discretized uniformly with time step size $\tau=0.05$ ms, resulting in a total of 20 time steps. 
We consider a longer time interval in Section \ref{sec: heartbeat} in order to test the robustness of the considered nonlinear solvers when simulating a complete heartbeat: in this case we simulate the electrical propagation for 0.5 second, thus performing 10'000 time steps. 

\begin{figure}[pb]
	\centering
	\includegraphics[scale=.2]{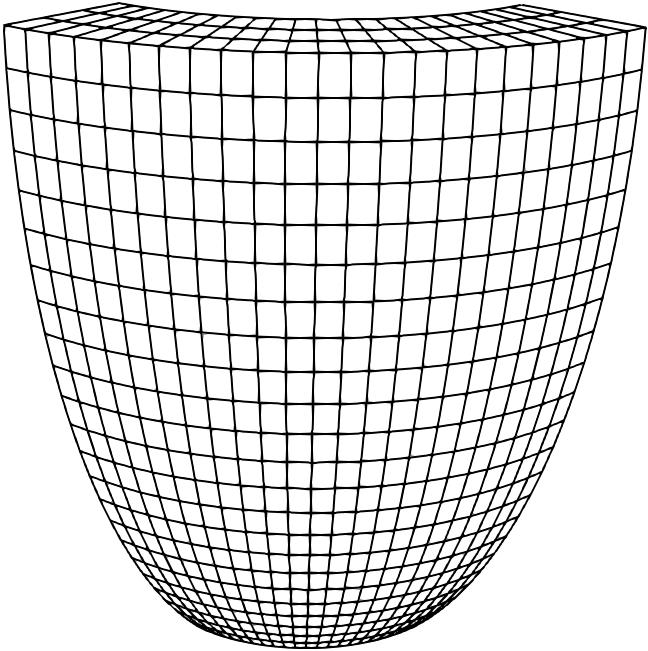}
	\caption{Idealized left ventricle geometry.}
	\label{fig: ellips}
\end{figure}

All nonlinear solvers are root-finding algorithms already implemented within the SNES (Nonlinear Algebraic Solvers) package contained in the PETSc (Portable, Extensible Toolkit for Scientific Computation) library \cite{balay2019petsc}. The tests are performed on the EOS supercomputer at the University of Pavia (\hyperref[https://matematica.unipv.it/cluster-di-calcolo/]{https://matematica.unipv.it/cluster-di-calcolo/}), composed of 672 Intel Xeon CPUs running at $2.1$ GHz.

\subsection{Tested nonlinear solvers.}
We do not restrict our numerical tests to the two nonlinear solvers investigated in the above theoretical study. The ones considered are the following: 
\begin{enumerate}
	\item the standard Newton method,  where the Jacobian system is solved with a preconditioned iterative method with a strict stopping criterion;
	\item inexact Newton, where the accuracy of the iterative solver increases with the iterations in order to avoid over-solving the tangent problem in the first iterations;
	\item quasi-Newton methods, where the Jacobian is either approximated by the action of a preconditioner or it is solved inexactly;
	\item nonlinear Generalized Minimal Residual (NGMRES) method, where an optimal mixing between the current iteration and the previous ones is computed at each iteration;
	\item nonlinear Conjugate Gradient (NCG) methods.
\end{enumerate}
Whenever preconditioning is needed, we consider the Algebraic Multigrid (AMG) implementation from PETSc library, recalled as PCGAMG from the PC (Preconditioners) package.
A sketch of the methods are given in the following paragraphs, together with the labels we will use in the results. We refer to the monograph \cite{wright1999numerical} for any further details.\\
In what follows, we consider for simplicity the arbitrary root-finding problem
\begin{equation}\label{eq:abstract equation}
	\vec F(\vec x) = \vec 0.
\end{equation} 

\subsubsection{Standard Newton \emph{[Newton-MG]}. }
As standard benchmark, we compare the performance of the above mentioned nonlinear solvers with the classic Newton method. It considers a first order approximation for the problem \eqref{eq:abstract equation}, which yields
\begin{equation}\label{eq:newton taylor}
	\vec F(\vec x^{k+1})\approx \vec F(\vec x^k) + \nabla_{\vec x}\vec F(\vec x^k)\Delta \vec x^{k+1}, 
	\qquad
	\Delta \vec x^{k+1} = \vec x^{k+1} - \vec x^k.
\end{equation}
By assuming that eventually $\vec F(\vec x^{k+1})=\vec 0$, we obtain the $k$-th tangent or Jacobian linear problem
\begin{equation}\label{eq:tangent newton problem}
	\nabla_{\vec x}\vec F(\vec x^k) \Delta \vec x^{k+1} = - \vec F(\vec x^k),
\end{equation}
which is usually solved iteratively with a strict stopping criterion (in our test, we use a multigrid (MG) method). Despite the quadratic convergence rate (which makes Newton method very appealing), this method requires to assemble and invert the gradient matrix $\nabla_{\vec x}$ which is a computationally expensive task.

\subsubsection{Inexact Newton with Eisenstat-Walker adaptive tolerance \emph{[iNewton]}. }
In order to overcome the computational costs required by the solution of the Jacobian system, it is possible to drop the accuracy in the first Newton iterations, resulting in lighter computational steps. Indeed, there is no need to solve accurately the approximation \eqref{eq:newton taylor} in the first iterations, since it yields a larger approximation error at the beginning of the iterative procedure. 
Moreover, this strategy can be further improved by considering some adaptive tolerances. For instance, the popular Eisenstat-Walker strategy \cite{eisenstat1994globally, eisenstat1996choosing} considers a decreasing tolerance that follows the validity of the first order approximation:
$$
\texttt{tol}^n = \frac{| \|\nabla_{\vec x}\vec F(\vec x^{k})[\Delta\vec x^k] - \vec F(\vec x^{k})\| - \|\vec F(\vec x^{k})\| |} {\|\vec F(\vec x^{k})\|},
$$
for a given initial tolerance $\texttt{tol}^0$. To our knowledge, this results in an overall increase of the nonlinear iterations with a significant reduction in the computational cost of each one of them. This algorithm is super-linearly convergent.

\subsubsection{Quasi-Newton methods \emph{[QN preonly, QN jac-low]}.}
This family of methods have already been introduced in Section \ref{sec: quasi-newton}. From a numerical point of view, we consider the limited-memory implementation, where only the last $m$ vectors ($\vec s$ and $\vec y$) are used, where $m$ is defined according to the tuning done in Section \ref{sec:tuning}. 
Additionally, as in \cite{barnafi2022parallel} we use the Jacobian approximation $\mat B_0$ to be equal to the exact Jacobian \rev{at each time step, i.e. $\mat B_0 = \grad_\vec x\vec F(\vec x^n)$}, and we consider the action of $[\mat B^0]^{-1}$ to be either solved approximately with only the action of the algebraic multigrid (AMG) preconditioner (preonly) or with 10 iterations of a CG iterative solver, preconditioned by an AMG (jac-low). We refer to Ref. \cite{barnafi2022parallel} for an alternative usage of inexact solvers, and thus a variant of the jac-low method, in the context of quasi-Newton methods by means of the relative tolerance of the linear solver instead of a fixed number of linear iterations. We highlight that both methods use require the assembly of the Jacobian at the current time step, i.e. the first Jacobian from a standard Netwon iteration.

\rev{
\begin{remark}
It is interesting to note that, at each timestep, we reassemble the Jacobian matrix, and thus consider the block $\int\partial_v I_\text{ion}(v, w)\phi_i\phi_j\,dx$ in it, which is not convex. Despite this, we have not observed lack of convergence because of this in any of our simulations. We believe this happens because the dynamics of this problem are very fast, which require small timesteps to obtain accurate solutions. This yields \eqref{eq:dt-bound} automatically, so there are no problems related to convexity in practice.
\end{remark}
}

\subsubsection{Nonlinear GMRES \emph{[NGMRES]}.}
This method arises as variant of the well-known linear GMRES. The underlying idea is simple: find, at each iteration, an optimal mixing between a new candidate given by a simple descent direction and the previous $m-1$ iterations \cite{washio1997krylov}.
This can be implemented by fixing a number $m$ of vectors to mix, then at iteration $k$ compute first a descent candidate $\vec x_k^M=\vec x_k - \vec F(\vec x_k)$ and then the mixing weights $\vec \alpha = (\alpha_0, ..., \alpha_{m-1})$, by minimizing the problem
$$ \min_{\vec \alpha} \left\|\vec F\left(\sum_{i=0}^{m-2}\alpha_i\vec x_{k-i} + \alpha_{m-1}\vec x_k^M\right) \right\|,
\qquad
\sum_{i=0}^{m-1}\alpha_i = 1.
$$ 
Additional implementation details can be found in Ref. \cite{brune2015composing}.

\subsubsection{Nonlinear Conjugate Gradient methods \emph{[NCG]}}\label{sec:ncg methods}
We refer to Section \ref{sec: nonlinear CG} and Ref. \cite{wright1999numerical} for a description of the Fletcher-Reeves (FR) method. The main NCG variants that we consider differ in the computation of the coefficient $\beta$ as follows: 
\begin{itemize}
	\item Fletcher--Reeves (FR): 
	$$ \beta_{k+1} = \frac{\grad f^T_{k+1}\cdot \grad f_{k+1}}{\grad f^T_{k}\cdot \grad f_{k}}.$$
	\item Polak-Ribière-Polyak (PRP):
	$$ \beta_{k+1} = \frac{\grad f^T_{k+1}\cdot(\grad f_{k+1} - \grad f_k)}{\grad f^T_{k}\cdot \grad f_{k}}.$$
	\item Dai–Yuan (DY):
	$$ \beta_{k+1} = \frac{\grad f^T_{k+1}\cdot \grad f_{k+1}}{-p_k(\grad f_{k+1} - \grad f_k)}.$$
	\item Conjugate-Descent (CD):
	$$ \beta_{k+1} = \frac{\grad f^T_{k+1}\cdot \grad f_{k+1}}{p_k\grad f_{k+1}}.$$
\end{itemize}


\subsection{Solver tuning}\label{sec:tuning}
In this section we proceed as in Ref. \cite{barnafi2022alternative} to obtain the parameters to be used in each method. We restrict the tuning to one parameter only to simplify the analysis, which we detail in what follows:
\begin{itemize}
	\item Inexact-Newton: we consider both the initial relative tolerance $\texttt{tol}^0$ to be 0.9, 0.5, 0.1 and 0.01.
	\item Quasi-Newton (QN) preonly and jac-low: we consider $m$ = 2, 5, 10 and 20 previous vectors for the low memory implementation.
	\item Nonlinear GMRES: we consider $m$ = 1,2,5 and 10 previous vectors.
	\item Nonlinear CG: we consider the updates described in Section \ref{sec:ncg methods}.
\end{itemize}

We consider a fixed number of $N_p=16$ processors and a fixed number of $64\times64\times64$ finite elements, resulting in roughly half million degrees of freedom (DoFs). Detailed instructions used for each method are listed in \ref{appendix:petsc}

The choice of parameters is reasonable but still largely arbitrary since there is still plenty of room for deeper studies of each method in order to obtain a truly optimal set of parameters. 
As a matter of fact, the choice of the geometry, as well as the ionic model or the healthy/pathological scenario can influence the choice of each method's parameters. 
This would require a fine-tuning of parameters for each different physical scenario, which is out of the scope of this work.

\begin{figure}
	\centering
	\begin{tikzpicture}
		\begin{axis}
			[width=0.98\textwidth, height=8cm, xlabel=Time $(ms)$, ylabel=CPU time (s), tick label style={font=\tiny},line width=1.5pt,legend style={draw=none, fill opacity=0.7, text opacity = 1,row sep=2pt, font=\small}, xmin=-0.05, xmax=2.05, legend pos=north east, legend columns=3, ymin=0, mark options={mark size=1.2pt}]
			\addplot+[\newtoncolor, mark=*] table [restrict x to domain=0:2, x=Time, y=SNEStime, col sep=comma, filter discard warning=false, unbounded coords=discard, each nth point=2] {results/heartbeat/HEARTBEAT_NEWTON_32.csv};
			\addplot+[\inewtoncolor, mark=square*]  table [restrict x to domain=0:2, x=Time, y=SNEStime, col sep=comma, filter discard warning=false, unbounded coords=discard, each nth point=2 ] {results/heartbeat/HEARTBEAT_NEWTON_EW_32.csv};
			\addplot+[\bfgscolor, mark=triangle*]  table [restrict x to domain=0:2, x=Time, y=SNEStime, col sep=comma, filter discard warning=false, unbounded coords=discard, each nth point=2 ] {results/heartbeat/HEARTBEAT_QN_PREONLY_32.csv};
			\addplot+[\ibfgscolor, mark=*]  table [restrict x to domain=0:2, x=Time, y=SNEStime, col sep=comma, filter discard warning=false, unbounded coords=discard, each nth point=2 ] {results/heartbeat/HEARTBEAT_QN_JACOBIAN_LOW_32.csv};
			\addplot+[\ngmrescolor, mark=square*]  table [restrict x to domain=0:2, x=Time, y=SNEStime, col sep=comma, filter discard warning=false, unbounded coords=discard, each nth point=2 ] {results/heartbeat/HEARTBEAT_NCG_32.csv};
			\addplot+[\ncgcolor, mark=triangle*]  table [restrict x to domain=0:2, x=Time, y=SNEStime, col sep=comma, filter discard warning=false, unbounded coords=discard, each nth point=2 ] {results/heartbeat/HEARTBEAT_NGMRES_32.csv};
			\legend {Newton-MG, iNewton, QN preonly, QN jac-low, NGMRES, NCG-PRP}
		\end{axis}
	\end{tikzpicture}
	\caption{\emph{Solver tuning, time evolution.} Fixed number of processors $N_p=16$ and fixed mesh of $64\times64\times64$ elements (550k DoFs). Ten Tusscher--Panfilov ionic model. Time evolution of global CPU times (in seconds).}
	\label{fig:sensitivity}
\end{figure}
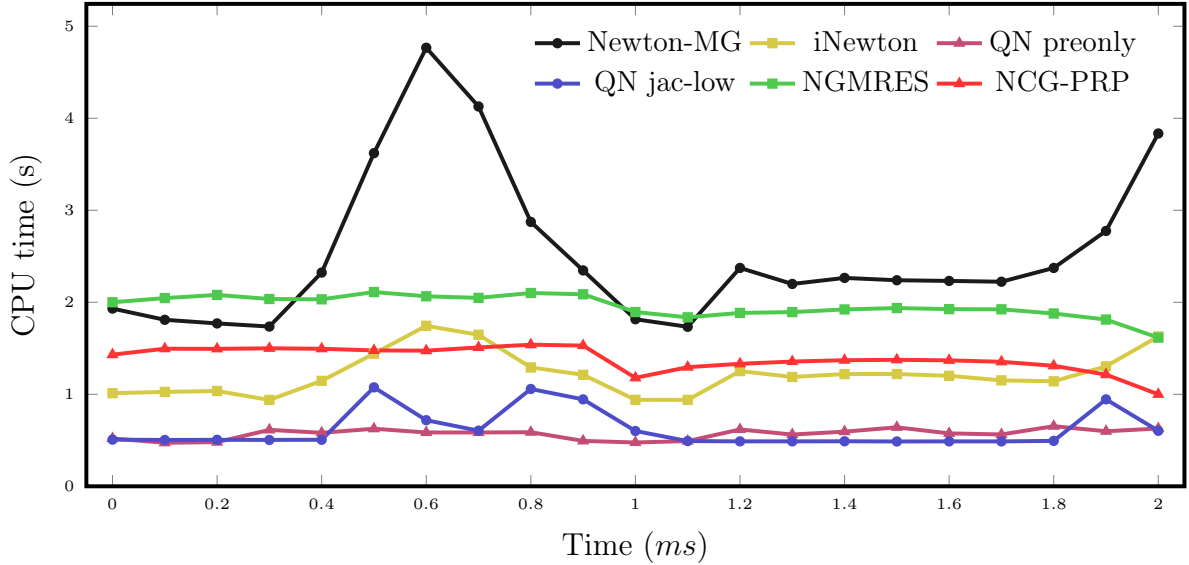

We consider the results obtained in Table \ref{tab:sensitivity}, where the best performing case is highlighted in bold fonts. We also report in Figure \ref{fig:sensitivity} a time evolution of the CPU times (in second). We keep these values as parameters for the next following tests, with only one exception: the QN preonly method. Indeed, the parameter $m$ yielding the best performance is too low ($m=2$), and this can cause convergence issues when the iterations required for convergence increase. For this reason, we will use $m=5$ vectors for simulations longer than a couple of milliseconds.
The following observations need to be stated:
\begin{itemize}
	\item As in Ref. \cite{barnafi2022alternative}, the inexact-Newton method performs best when an adequate initial tolerance is found, so that the number of nonlinear iterations is very similar to standard Newton. Indeed, we found that the optimal tolerance is the largest one such that the average nonlinear iterations are equal to those incurred by the Newton-MG method.
	\item A distinctive characteristic of the first order methods (NCG, NGMRES) is that they require much more iterations than the other methods to converge.
	\item Common choices of nonlinear solvers for the Bidomain equations involve linearization, i.e. IMEX schemes or one Newton iteration. If we consider the latter for comparison, we highlight that Quasi-Newton methods are vastly superior, as the times it takes to do one Newton iteration ($620.75/5.05\approx 120$ seconds) is comparable to the time it takes to solve the entire nonlinear problem ($122.94$ seconds with QN jac-low, $86.10$ seconds with QN preonly). 
	\item First order methods are faster than standard Newton, but still they are not competitive against the other methods under consideration. Their attractiveness resides in the fact that they do not require the assembly of a Jacobian, making them very well-suited for matrix-free frameworks, in particular using GPU accelerated computing.
\end{itemize}

\begin{table}[pb]
	\caption{\emph{Solver tuning.} Average number of nonlinear iterations and global CPU times (in seconds) for the nonlinear Bidomain solvers considered. The best performance is highlighted in bold font. Fixed number of processors $N_p=16$ and fixed mesh of $64\times64\times64$ elements (550k DoFs). Ten Tusscher--Panfilov ionic model. }
	\label{tab:sensitivity}
	\centering
	\begin{tabular}{@{}lrr@{}}
		\toprule Method & Average iterations & CPU time (s) \\
		\midrule {\bf Newton-MG} & \textbf{5.05} & \textbf{620.75} \\
		\midrule 
		iNewton (rtol=$0.001$) & 5.05 & 528.56 \\
		iNewton (rtol=$0.01$) & 5.05 &  225.97 \\
		\textbf{iNewton (rtol=$0.1$)} & \textbf{5.05} &  \textbf{208.06} \\
		iNewton (rtol=$0.5$) & 6.1 & 230.94 \\
		\midrule 
		\textbf{QN preonly $(m=2)$}  & \textbf{27.70} & \textbf{86.10} \\
		QN preonly $(m=5)$  & 31.70 & 92.78 \\
		QN preonly $(m=10)$ & 34.20 & 98.83 \\
		QN preonly $(m=20)$ & 38.25 & 113.18 \\
		\midrule 
		QN jac-low ($m=2$)  & 5.85 & 148.65 \\
		QN jac-low ($m=5$)  & 6.55 & 168.90 \\
		QN jac-low ($m=10$) & 5.10 & 127.00 \\
		\textbf{QN jac-low ($m=20$)}& \textbf{4.6} & \textbf{122.94} \\ 
		\midrule 
		\textbf{NGMRES ($m=2$)}   & \textbf{545.75} & \textbf{281.09} \\
		NGMRES ($m=5) $  & 545.75 & 298.63 \\
		NGMRES ($m=10)$  & 545.75 & 317.33 \\
		NGMRES ($m=20)$  & 545.75 & 340.23 \\
		\midrule 
		NCG (FR)  & 917.35 & 486.67 \\
		\textbf{NCG (PRP)} & \textbf{917.3} & \textbf{448.95} \\
		NCG (DY)  & 917.3 & 454.37 \\
		NCG (CD)  & 917.3 & 464.83 \\\bottomrule
	\end{tabular}
\end{table}

\subsection{Robustness with respect to the problem size}
In this section, we study the robustness of the solvers under consideration with respect to the number of degrees of freedom. We consider a time interval of $[0,1]$ ms and the Ten Tusscher--Panfilov ionic model. The number of processors is fixed to $N_p=16$ and the mesh size increases from $32\times32\times32$ to $102\times102\times102$ elements, i.e. from 72k to 2M DoFs. The results are shown in Figures \ref{fig:robustness newton inewton}, \ref{fig:robustness ngmres ncg} and \ref{fig:robustness qn} and the average CPU times and iterations are collected in Table \ref{tab:robustness dofs}. The methods have been grouped according to their behavior: both Newton methods present a mild increase of iterations towards the end of the first millisecond, whereas at the beginning they present a more constant behavior. First order methods (NGMRES, NCG-PRP) are not robust with respect to the problem size, and the dependency is much more severe for the NCG-PRP, which presents an increase of over 200 iterations between the coarsest and finest cases. Finally, quasi-Newton methods are the most robust ones, and indeed there is no appreciable trend in the iteration counts. 

In view of this, we conclude that all (quasi, inexact and standard) Newton methods are robust with respect to the problem size, with the most standard approaches presenting only a variation of one or two iterations, while first order methods are less robust.

\def\figHeight{6cm}
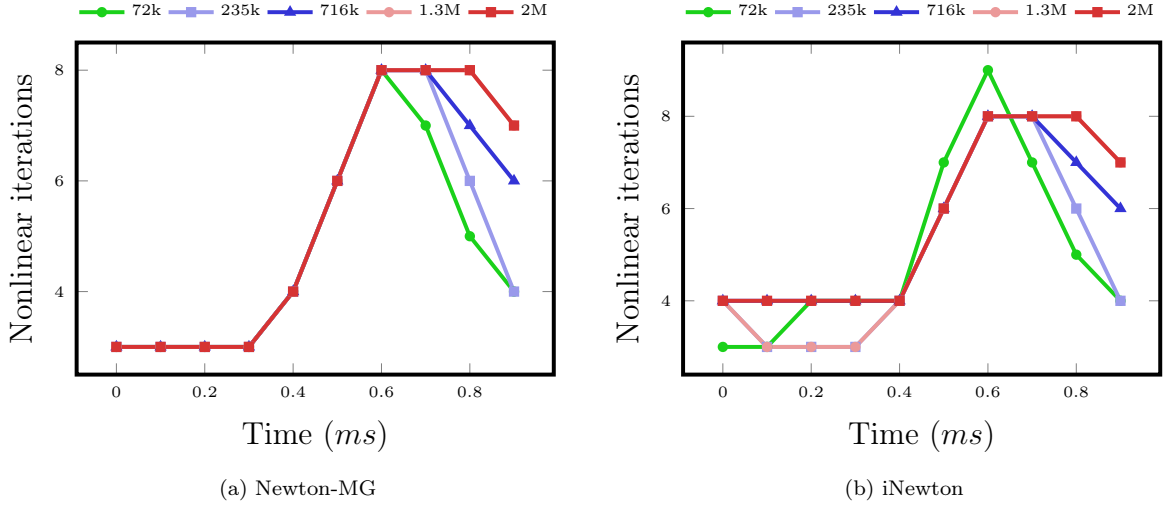
\begin{figure}[pb]
	\centering
    \begin{subfigure}[b]{0.48\textwidth}
	    \begin{tikzpicture}
	    	\begin{axis}[width=\linewidth, height=\figHeight, xlabel=Time $(ms)$, ylabel=Nonlinear iterations, tick label style={font=\tiny},line width=1.5pt, legend style={draw=none, fill opacity=0.7, text opacity = 1,row sep=2pt,font=\tiny, at={(0.5,1.03)}, anchor=south}, legend columns=5, mark options={mark size=1.2pt}]
	    		\addplot+[green!80!black!90!white, mark=*]  table [x=Time, y=snesIts, col sep=comma, each nth point=2] {results/robustness-size/ROBUSTNESS_NEWTON_32.csv};
	    		\addplot+[blue!80!black!40!white, mark=square*]  table [x=Time, y=snesIts, col sep=comma, each nth point=2] {results/robustness-size/ROBUSTNESS_NEWTON_48.csv};
	    		\addplot+[blue!80!black!80!white, mark=triangle*]  table [x=Time, y=snesIts, col sep=comma, each nth point=2] {results/robustness-size/ROBUSTNESS_NEWTON_70.csv};
	    		\addplot+[red!80!black!40!white, mark=*]  table [x=Time, y=snesIts, col sep=comma, each nth point=2] {results/robustness-size/ROBUSTNESS_NEWTON_86.csv};
	    		\addplot+[red!80!black!80!white, mark=square*]  table [x=Time, y=snesIts, col sep=comma, each nth point=2] {results/robustness-size/ROBUSTNESS_NEWTON_102.csv};
	    		\legend {72k,235k,716k,1.3M,2M}
	    	\end{axis}
	    \end{tikzpicture}
        \caption{Newton-MG}
    \end{subfigure}
    \begin{subfigure}[b]{0.48\textwidth}
	    \begin{tikzpicture}
	    	\begin{axis} [width=\linewidth, height=\figHeight, xlabel=Time $(ms)$, ylabel=Nonlinear iterations, tick label style={font=\tiny},line width=1.5pt, legend style={draw=none, fill opacity=0.7, text opacity = 1,row sep=2pt,font=\tiny, at={(0.5,1.03)}, anchor=south}, legend columns=5, mark options={mark size=1.2pt}]
	    		\addplot+[green!80!black!90!white, mark=*] table [x=Time, y=snesIts, col sep=comma, each nth point=2] {results/robustness-size/ROBUSTNESS_NEWTON_EW_32.csv};
	    		\addplot+[blue!80!black!40!white, mark=square*]  table [x=Time, y=snesIts, col sep=comma, each nth point=2] {results/robustness-size/ROBUSTNESS_NEWTON_EW_48.csv};
	    		\addplot+[blue!80!black!80!white, mark=triangle*] table [x=Time, y=snesIts, col sep=comma, each nth point=2] {results/robustness-size/ROBUSTNESS_NEWTON_EW_70.csv};
	    		\addplot+[red!80!black!40!white, mark=*]  table [x=Time, y=snesIts, col sep=comma, each nth point=2] {results/robustness-size/ROBUSTNESS_NEWTON_EW_86.csv};
	    		\addplot+[red!80!black!80!white, mark=square*] table [x=Time, y=snesIts, col sep=comma, each nth point=2] {results/robustness-size/ROBUSTNESS_NEWTON_EW_102.csv};
	    		\legend {72k,235k,716k,1.3M,2M}
	    	\end{axis}
	    \end{tikzpicture}
        \caption{iNewton}
    \end{subfigure}
	\caption{\emph{Robustness of (a) Newton-MG and (b) inexact-Newton solvers with respect to the problem size.} Fixed number of processors $N_p=16$ and increasing number of degrees of freedom from 72k to 2M. Nonlinear iterations over the time interval $[0,1]$ ms.}
	\label{fig:robustness newton inewton}
\end{figure}

\begin{figure}[pb]
	\centering
    \begin{subfigure}[b]{0.48\textwidth}
	    \begin{tikzpicture}
	    	\begin{axis}[width=\linewidth, height=\figHeight, xlabel=Time $(ms)$, ylabel=Nonlinear iterations, tick label style={font=\tiny},line width=1.5pt, legend style={draw=none, fill opacity=0.7, text opacity = 1,row sep=2pt,font=\tiny, at={(0.5,1.03)}, anchor=south}, legend columns=5, mark options={mark size=1.2pt}]
	    		\addplot+[green!80!black!90!white, mark=*]  table [x=Time, y=snesIts, col sep=comma, each nth point=2] {results/robustness-size/ROBUSTNESS_QN_PREONLY_32.csv};
	    		\addplot+[blue!80!black!40!white, mark=square*]  table [x=Time, y=snesIts, col sep=comma, each nth point=2] {results/robustness-size/ROBUSTNESS_QN_PREONLY_48.csv};
	    		\addplot+[blue!80!black!80!white, mark=triangle*]  table [x=Time, y=snesIts, col sep=comma, each nth point=2] {results/robustness-size/ROBUSTNESS_QN_PREONLY_70.csv};
	    		\addplot+[red!80!black!40!white, mark=*]  table [x=Time, y=snesIts, col sep=comma, each nth point=2] {results/robustness-size/ROBUSTNESS_QN_PREONLY_86.csv};
	    		\addplot+[red!80!black!80!white, mark=square*]  table [x=Time, y=snesIts, col sep=comma, each nth point=2] {results/robustness-size/ROBUSTNESS_QN_PREONLY_102.csv};
	    		\legend {72k,235k, 716k, 1.3M, 2M}
	    	\end{axis}
	    \end{tikzpicture}
        \caption{QN preonly}
    \end{subfigure}
    \begin{subfigure}[b]{0.48\textwidth}
	    \begin{tikzpicture}
	    	\begin{axis}[width=\linewidth, height=\figHeight, xlabel=Time $(ms)$, ylabel=Nonlinear iterations, tick label style={font=\tiny},line width=1.5pt, legend style={draw=none, fill opacity=0.7, text opacity = 1,row sep=2pt,font=\tiny, at={(0.5,1.03)}, anchor=south}, legend columns=5, mark options={mark size=1.2pt}]
	    		\addplot+ [green!80!black!90!white, mark=*]  table [x=Time, y=snesIts, col sep=comma, each nth point=2] {results/robustness-size/ROBUSTNESS_QN_JACOBIAN_LOW_32.csv};
	    		\addplot+ [blue!80!black!40!white, mark=square*]  table [x=Time, y=snesIts, col sep=comma, each nth point=2] {results/robustness-size/ROBUSTNESS_QN_JACOBIAN_LOW_48.csv};
	    		\addplot+ [blue!80!black!80!white, mark=triangle*]  table [x=Time, y=snesIts, col sep=comma, each nth point=2] {results/robustness-size/ROBUSTNESS_QN_JACOBIAN_LOW_70.csv};
	    		\addplot+ [red!80!black!40!white, mark=*]  table [x=Time, y=snesIts, col sep=comma, each nth point=2] {results/robustness-size/ROBUSTNESS_QN_JACOBIAN_LOW_86.csv};
	    		\addplot+ [red!80!black!80!white, mark=square*]  table [x=Time, y=snesIts, col sep=comma, each nth point=2] {results/robustness-size/ROBUSTNESS_QN_JACOBIAN_LOW_102.csv};
	    		\legend {72k,235k, 716k, 1.3M, 2M}
	    	\end{axis}
	    \end{tikzpicture}
        \caption{QN jac-low}
    \end{subfigure}
	\caption{\emph{Robustness of (a) QN preonly  and (b) QN jac-low solvers with respect to the problem size.} Fixed number of processors $N_p=16$ and increasing number of degrees of freedom from 72k to 2M. Nonlinear iterations over the time interval $[0,1]$ ms.}
	\label{fig:robustness qn}
\end{figure}
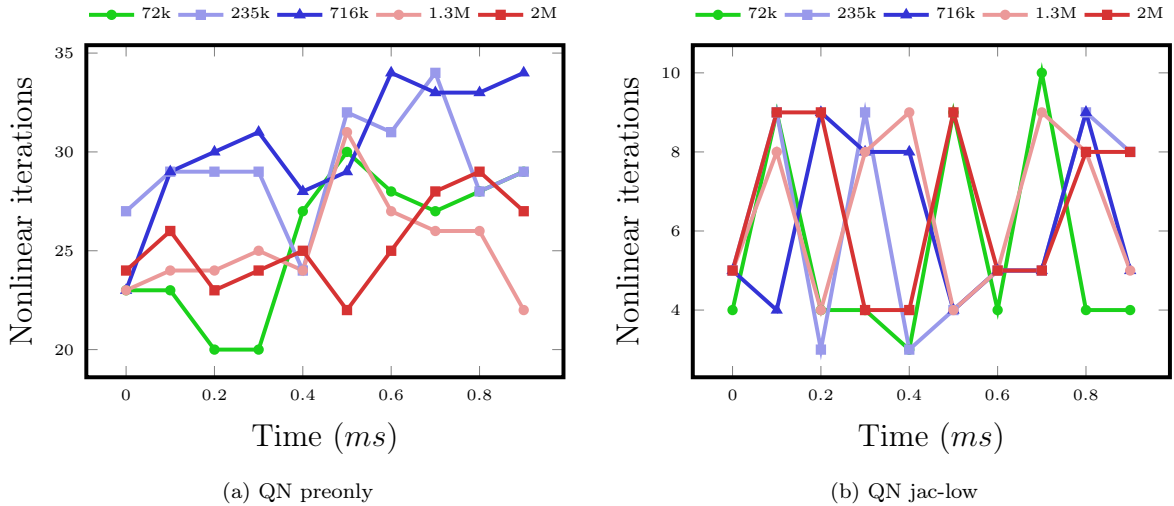

\begin{figure}[pb]
	\centering
    \begin{subfigure}[b]{0.48\textwidth}
	    \begin{tikzpicture}
	    	\begin{axis}[width=\linewidth, height=\figHeight, xlabel=Time $(ms)$, ylabel=Nonlinear iterations, tick label style={font=\tiny},line width=1.5pt, legend style={draw=none, fill opacity=0.7, text opacity = 1,row sep=2pt,font=\tiny, at={(0.5,1.03)}, anchor=south}, legend columns=5, mark options={mark size=1.2pt}]
	    		\addplot+ [green!80!black!90!white, mark=*] table [x=Time, y=snesIts, col sep=comma, each nth point=2] {results/robustness-size/ROBUSTNESS_NGMRES_32.csv};
	    		\addplot+ [blue!80!black!40!white, mark=square*] table [x=Time, y=snesIts, col sep=comma, each nth point=2] {results/robustness-size/ROBUSTNESS_NGMRES_48.csv};
	    		\addplot+ [blue!80!black!80!white, mark=triangle*] table [x=Time, y=snesIts, col sep=comma, each nth point=2] {results/robustness-size/ROBUSTNESS_NGMRES_70.csv};
	    		\addplot+ [red!80!black!40!white, mark=*]  table [x=Time, y=snesIts, col sep=comma, each nth point=2] {results/robustness-size/ROBUSTNESS_NGMRES_86.csv};
	    		\addplot+ [red!80!black!80!white, mark=square*] table [x=Time, y=snesIts, col sep=comma, each nth point=2] {results/robustness-size/ROBUSTNESS_NGMRES_102.csv};
	    		\legend {72k,235k, 716k, 1.3M, 2M}
	    	\end{axis}
	    \end{tikzpicture}
        \caption{NGMRES}
    \end{subfigure}
    \begin{subfigure}[b]{0.48\textwidth}
	    \begin{tikzpicture}
	    	\begin{axis}
	    		[width=\linewidth, height=\figHeight, xlabel=Time $(ms)$, ylabel=Nonlinear iterations, tick label style={font=\tiny},line width=1.5pt, legend style={draw=none, fill opacity=0.7, text opacity = 1,row sep=2pt,font=\tiny, at={(0.5,1.03)}, anchor=south}, legend columns=5, mark options={mark size=1.2pt}]
	    		\addplot+ [green!80!black!90!white, mark=*] table [x=Time, y=snesIts, col sep=comma, each nth point=2] {results/robustness-size/ROBUSTNESS_NCG_32.csv};
	    		\addplot+ [blue!80!black!40!white, mark=square*] table [x=Time, y=snesIts, col sep=comma, each nth point=2] {results/robustness-size/ROBUSTNESS_NCG_48.csv};
	    		\addplot+ [blue!80!black!80!white, mark=triangle*] table [x=Time, y=snesIts, col sep=comma, each nth point=2] {results/robustness-size/ROBUSTNESS_NCG_70.csv};
	    		\addplot+ [red!80!black!40!white, mark=*]  table [x=Time, y=snesIts, col sep=comma, each nth point=2] {results/robustness-size/ROBUSTNESS_NCG_86.csv};
	    		\addplot+ [red!80!black!80!white, mark=square*] table [x=Time, y=snesIts, col sep=comma, each nth point=2] {results/robustness-size/ROBUSTNESS_NCG_102.csv};
	    		\legend {72k,235k, 716k, 1.3M, 2M}
	    	\end{axis}
	    \end{tikzpicture}
        \caption{NCG-PRP}
    \end{subfigure}
	\caption{\emph{Robustness of NGMRES (left) NCG-PRP (right) solvers with respect to the problem size.} Fixed number of processors $N_p=16$ and increasing number of degrees of freedom from 72k to 2M. Nonlinear iterations over the time interval $[0,1]$ ms.}
	\label{fig:robustness ngmres ncg}
\end{figure}
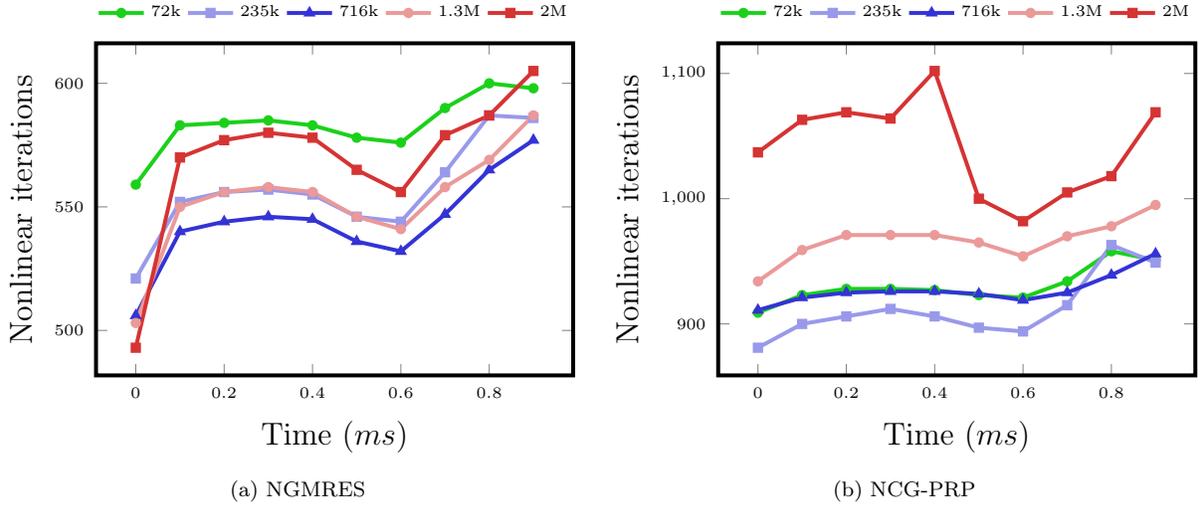

\begin{table}[pb]
	\caption{\emph{Nonlinear Bidomain solvers robustness with respect to the problem size.} Fixed number of processors $N_p=16$ and increasing number of degrees of freedom (DoFs) from 72k to 2M. Total CPU times in seconds (time) and average nonlinear iterations (nit) over the time interval $[0,1]$ ms. }
	\label{tab:robustness dofs}
	\centering
	\small
	\begin{tabular}{@{}rrrrrrrrrrr@{}}
		\toprule
		\multirow{2}{*}{Solver}		&\multicolumn{2}{c}{72k} &  \multicolumn{2}{c}{235k} &\multicolumn{2}{c}{716k}&\multicolumn{2}{c}{1.3M}&\multicolumn{2}{c}{2M}	\\
		\cmidrule{2-11}
		&time&nit&time&nit&time&nit&time&nit&time&nit\\
		\midrule
		Newton-MG  & 44.0 & 4.6 & 211.6 & 4.9 & 741.9 & 5.2 & 1391.0 & 5.4 & 2310.8 & 5.4\\
		iNewton    & 22.5 & 4.7 & 88.3 & 4.9 & 289.2 & 5.4 & 518.7 & 5.4 & 920.6 & 5.7 \\
		QN preonly & 8.1 & 25.8 & 33.2 & 29.4 & 116.6 & 29.8 & 204.9 & 27.3 & 327.6 & 25.4 \\
		QN jac-low & 11.5 & 4.2 & 59.6 & 5.8 & 142.2 & 4.1 & 335.1 & 5.3 & 561.8 & 5.3 \\
		NGMRES     & 24.5 & 584.9 & 113.1 & 559.0 & 376.3 & 546.0 & 774.0 & 555.0 & 1421.2 & 572.6 \\
		NCG-PRP    & 34.6 & 930.7 & 180.4 & 913.4 & 626.9 & 928.4 & 1208.3 & 968.2 & 2205.6 & 1040.6 \\
		\bottomrule
	\end{tabular}
\end{table}

\subsection{Impact of localized conductivity reduction (ischemia)}
\begin{figure}[pb]
	\centering
    \begin{subfigure}{0.49\textwidth}
	    \begin{tikzpicture}
	    	\begin{axis}
	    		[width=\linewidth, height=\figHeight, xlabel=Time $(ms)$, ylabel=Nonlinear iterations, tick label style={font=\tiny},line width=1.5pt, legend style={draw=none, fill opacity=0.7, text opacity = 1,row sep=2pt,font=\tiny, at={(0.5,1.03)}, anchor=south}, legend columns=5, mark options={mark size=1.2pt}]
	    		\addplot+[mark=*]  table [x=Time, y=snesIts, col sep=comma] {results/robustness-size/ROBUSTNESS_NGMRES_102.csv};
	    		\addplot+[dashed]  table [x=Time, y=snesIts, col sep=comma] {results/robustness-ischemic/ROBUSTNESS_NGMRES_102.csv};
	    		\legend {healthy, ischemic}
	    	\end{axis}
	    \end{tikzpicture}
        \caption{NGMRES}
    \end{subfigure}
    \begin{subfigure}{0.49\textwidth}
	    \begin{tikzpicture}
	    	\begin{axis}
	    		[width=\linewidth, height=\figHeight, xlabel=Time $(ms)$, ylabel=Nonlinear iterations, tick label style={font=\tiny},line width=1.5pt, legend style={draw=none, fill opacity=0.7, text opacity = 1,row sep=2pt,font=\tiny, at={(0.5,1.03)}, anchor=south}, legend columns=5, mark options={mark size=1.2pt}]
	    		\addplot+[mark=*]  table [x=Time, y=snesIts, col sep=comma] {results/robustness-size/ROBUSTNESS_NCG_102.csv};
	    		\addplot+[dashed]  table [x=Time, y=snesIts, col sep=comma] {results/robustness-ischemic/ROBUSTNESS_NCG_102.csv};
	    		\legend {healthy, ischemic}
	    	\end{axis}
	    \end{tikzpicture}
        \caption{NCG}
    \end{subfigure}
	\caption{\emph{Robustness of NGMRES (left) and NCG-PRP (right) solvers with respect to discontinuities in conductivity (ischemic scenario).} Fixed number of processors $N_p=16$ and approx 2 millions degrees of freedom. Nonlinear iterations over the time interval $[0,1]$ ms.}
	\label{fig:robustness ischemic ngmres ncg}
\end{figure}
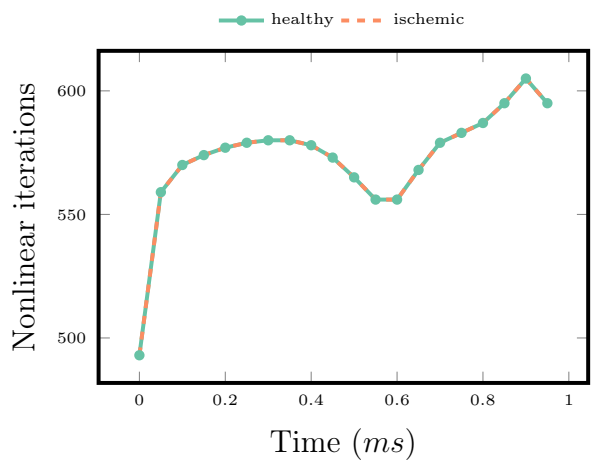
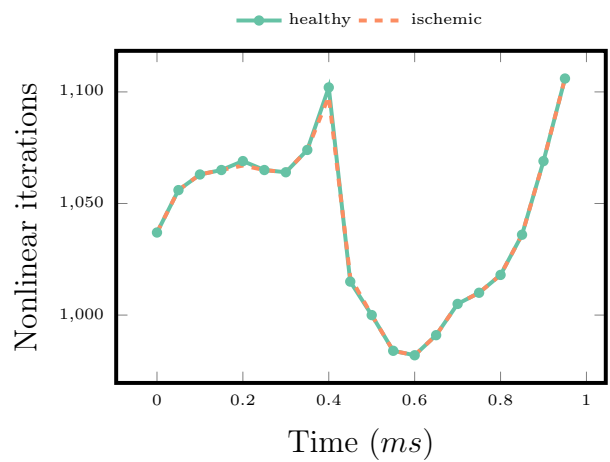

In this section, we compare the performance of the solvers in presence of an ischemic region in the cardiac tissue. The small, regular ischemic region is intramural, i.e. it runs from epi- to endocardium, is positioned in the middle of the considered geometry, and it presents reduced conductivity coefficients, see Table \ref{tab: physio coeff ischemic}. Moreover, the potassium extracellular concentration $K_o$ is increased from 5.4 mV to 8 mV, and the sodium conductance $G_{Na}$ is decreased by $30\%$, simulating a region with moderate ischemia. 
The number of processors is fixed to $N_p=16$ and the mesh consists of $102\times102\times102$ elements, resulting in approximatively 2M DoFs. We show the results in Figures \ref{fig:robustness ischemic newton}, \ref{fig:robustness ischemic qn} and \ref{fig:robustness ischemic ngmres ncg}.

We first highlight that Newton-MG, inexact-Newton and first order (NGMRES, NCG-PRP) methods are robust with respect to the presence of an ischemic region, as their nonlinear iterations do not change. This might not be so surprising for the Newton methods, whose number of linear iterations change as lower diffusions yield worse conditioning of the tangent system, but instead for first order methods this result is very interesting. Indeed, this can be inferred from the convergence theory developed in Section \ref{sec: convergence nonlinear CG}, where the conductivities play no role in the estimates used for establishing the convergence of the method. On the other hand, this result is not true for quasi-Newton methods, since QN preonly presents a performance deterioration throughout the entire time lapse considered. We note that the QN convexity estimate \eqref{eq:dt-bound} originates from an inequality where the conductivities provide a contribution, so that in fact a performance deterioration of the quasi-Newton methods can be expected, as the overall convexity of the potential decreases.

\begin{table}[ht]
	\caption{Conductivity coefficients for the Bidomain model in physiological and ischemic tissue.}
	\label{tab: physio coeff ischemic}
	\centering
	\begin{tabular}{@{}rcc@{}}
		\toprule
		Test						& $\sigma_l^i$		& $\sigma_t^i$			\\
		\midrule
		Normal 						& $3 \times 10^{-3} \Omega^{-1} \text{ cm}^{-1}$ 		& $3.1525 \times 10^{-4} \Omega^{-1} \text{ cm}^{-1}$ 	\\
		Ischemic 		& $1.5 \times 10^{-3} \Omega^{-1} \text{ cm}^{-1}$	   & $1.57625 \times 10^{-4} \Omega^{-1} \text{ cm}^{-1}$  \\
		\bottomrule
	\end{tabular}
\end{table}

\begin{figure}[pb]
	\centering
    \begin{subfigure}{0.49\textwidth}
	    \begin{tikzpicture}
	    	\begin{axis}
	    		[width=\linewidth, height=\figHeight, xlabel=Time $(ms)$, ylabel=Nonlinear iterations, tick label style={font=\tiny},line width=1.5pt, legend style={draw=none, fill opacity=0.7, text opacity = 1,row sep=2pt,font=\tiny, at={(0.5,1.03)}, anchor=south}, legend columns=5, mark options={mark size=1.2pt}]
	    		\addplot+[mark=*]       table [x=Time, y=snesIts, col sep=comma] {results/robustness-size/ROBUSTNESS_NEWTON_102.csv};
	    		\addplot+[dashed]  table [x=Time, y=snesIts, col sep=comma] {results/robustness-ischemic/ROBUSTNESS_NEWTON_102.csv};
	    		\legend {healthy, ischemic}
	    	\end{axis}
	    \end{tikzpicture}
        \caption{Newton}
    \end{subfigure}
    \begin{subfigure}{0.49\textwidth}
	    \begin{tikzpicture}
	    	\begin{axis}
	    		[width=\linewidth, height=\figHeight, xlabel=Time $(ms)$, ylabel=Nonlinear iterations, tick label style={font=\tiny},line width=1.5pt, legend style={draw=none, fill opacity=0.7, text opacity = 1,row sep=2pt,font=\tiny, at={(0.5,1.03)}, anchor=south}, legend columns=5, mark options={mark size=1.2pt}]
	    		\addplot+[mark=*]  table [x=Time, y=snesIts, col sep=comma] {results/robustness-size/ROBUSTNESS_NEWTON_EW_102.csv};
	    		\addplot+[dashed]  table [x=Time, y=snesIts, col sep=comma] {results/robustness-ischemic/ROBUSTNESS_NEWTON_EW_102.csv};
	    		\legend {healthy, ischemic}
	    	\end{axis}
	    \end{tikzpicture}
        \caption{Newton}
    \end{subfigure}
	\caption{\emph{Robustness of Newton-MG (left) and iNewton (right) solvers with respect to discontinuities in conductivity (ischemic scenario).} Fixed number of processors $N_p=16$ and approx 2 millions degrees of freedom. Nonlinear iterations over the time interval $[0,1]$ ms.}
	\label{fig:robustness ischemic newton}
\end{figure}

\begin{figure}[pb]
	\centering
    \begin{subfigure}{0.49\textwidth}
	    \begin{tikzpicture}
	    	\begin{axis}
	    		[width=\linewidth, height=\figHeight, xlabel=Time $(ms)$, ylabel=Nonlinear iterations, tick label style={font=\tiny},line width=1.5pt, legend style={draw=none, fill opacity=0.7, text opacity = 1,row sep=2pt,font=\tiny, at={(0.5,1.03)}, anchor=south}, legend columns=5, mark options={mark size=1.2pt}]
	    		\addplot+[mark=*]  table [x=Time, y=snesIts, col sep=comma] {results/robustness-size/ROBUSTNESS_QN_PREONLY_102.csv};
	    		\addplot+[dashed]  table [x=Time, y=snesIts, col sep=comma] {results/robustness-ischemic/ROBUSTNESS_QN_PREONLY_102.csv};
	    		\legend {healthy, ischemic}
	    	\end{axis}
	    \end{tikzpicture}
        \caption{QN-preonly}
    \end{subfigure}
    \begin{subfigure}{0.49\textwidth}
	    \begin{tikzpicture}
	    	\begin{axis}
	    		[width=\linewidth, height=\figHeight, xlabel=Time $(ms)$, ylabel=Nonlinear iterations, tick label style={font=\tiny},line width=1.5pt, legend style={draw=none, fill opacity=0.7, text opacity = 1,row sep=2pt,font=\tiny, at={(0.5,1.03)}, anchor=south}, legend columns=5, mark options={mark size=1.2pt}]
	    		\addplot+[mark=*]  table [x=Time, y=snesIts, col sep=comma] {results/robustness-size/ROBUSTNESS_QN_JACOBIAN_LOW_102.csv};
	    		\addplot+[dashed]  table [x=Time, y=snesIts, col sep=comma] {results/robustness-ischemic/ROBUSTNESS_QN_JACOBIAN_LOW_102.csv};
	    		\legend {healthy, ischemic}
	    	\end{axis}
	    \end{tikzpicture}
        \caption{QN-jac low}
    \end{subfigure}
	\caption{\emph{Robustness of QN preonly (left) and QN jac-low (right) solvers with respect to discontinuities in conductivity (ischemic scenario).} Fixed number of processors $N_p=16$ and approx 2 millions degrees of freedom. Nonlinear iterations over the time interval $[0,1]$ ms.}
	\label{fig:robustness ischemic qn}
\end{figure}
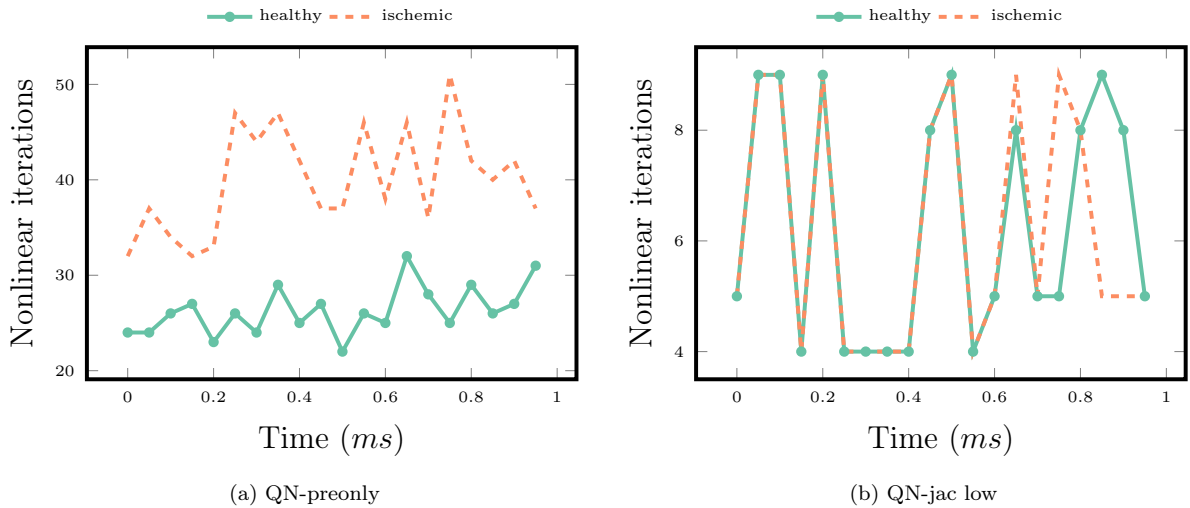

\subsection{Full activation-recovery simulation} \label{sec: heartbeat}
\begin{figure}[pb]
	\begin{center}
		\begin{multicols}{3}
			t = 1 ms 	\vspace*{2mm} \\	\includegraphics[scale=.18]{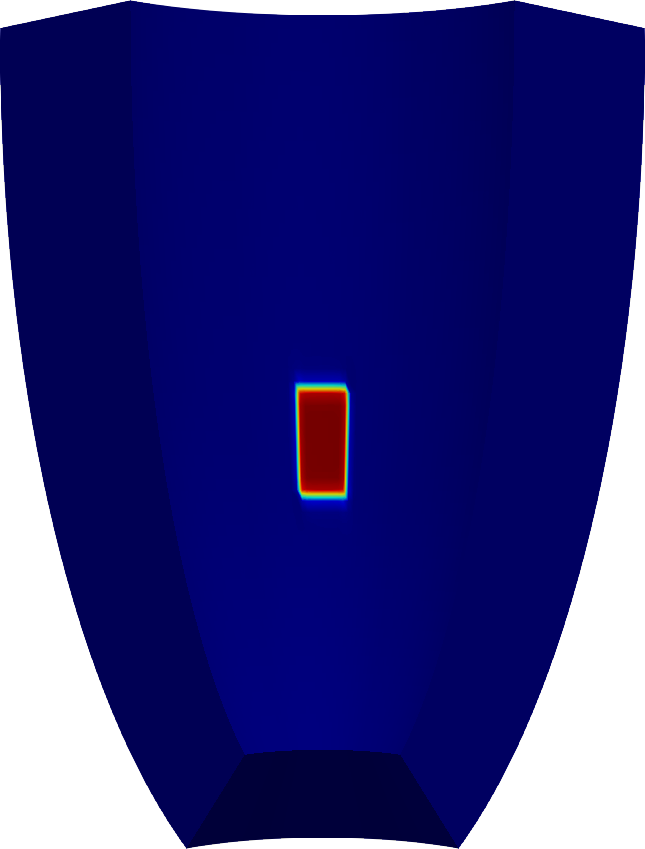} \\	\vspace*{2mm}
			\columnbreak
			t = 20 ms 	\vspace*{2mm} \\	\includegraphics[scale=.18]{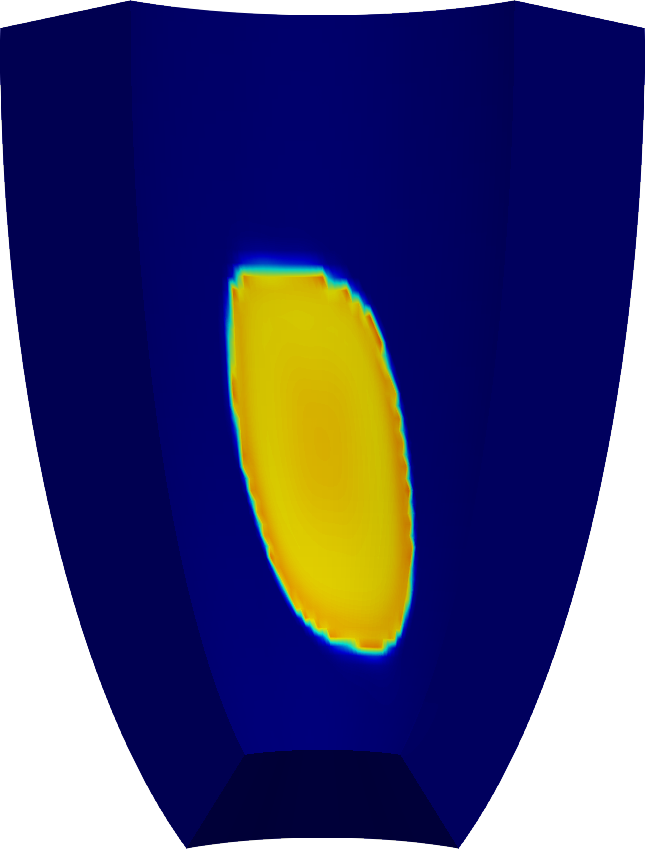} \\	\vspace*{2mm}
			\columnbreak
			t = 40 ms 	\vspace*{2mm} \\	\includegraphics[scale=.18]{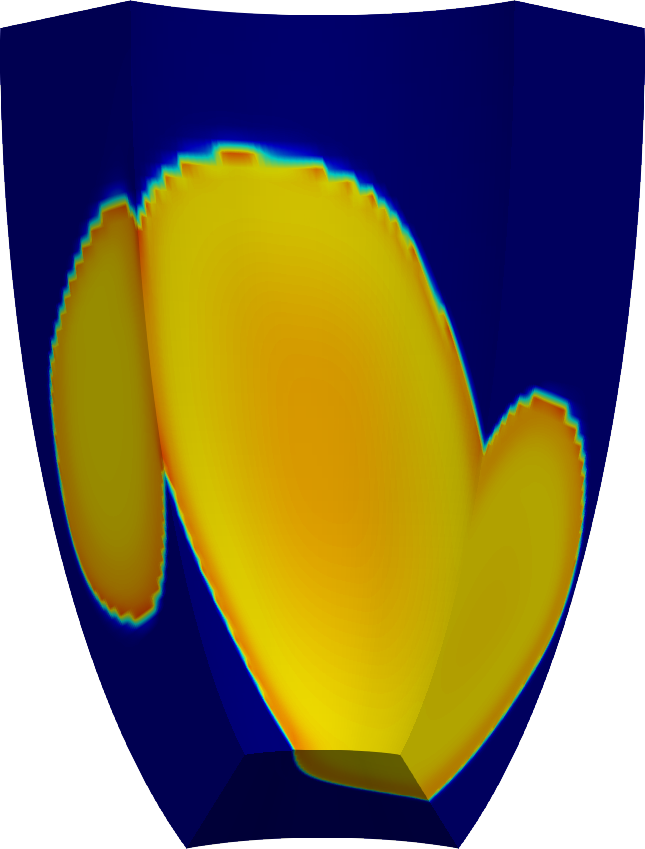} \\	\vspace*{2mm}
		\end{multicols}
		\includegraphics[scale=.40]{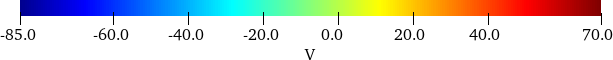}
		\begin{multicols}{3}
			t = 1 ms 	\vspace*{2mm} \\	\includegraphics[scale=.18]{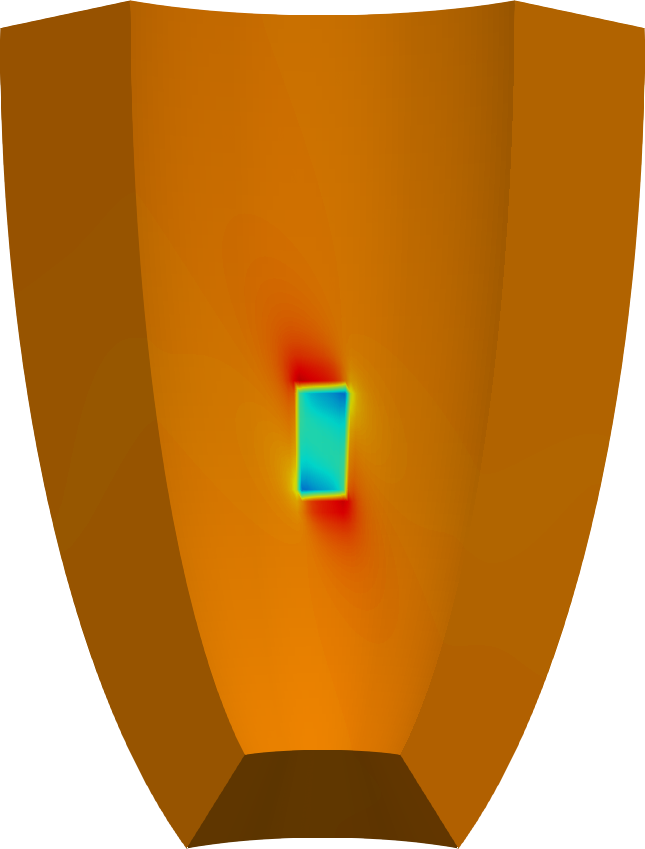} \\	\vspace*{2mm}
			\columnbreak
			t = 20 ms 	\vspace*{2mm} \\	\includegraphics[scale=.18]{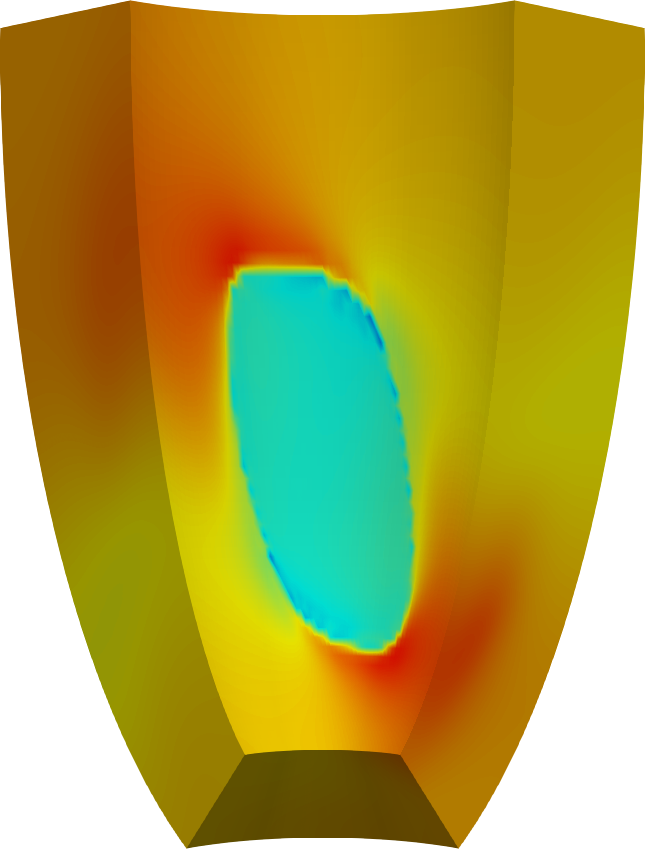} \\	\vspace*{2mm}
			\columnbreak
			t = 40 ms 	\vspace*{2mm} \\	\includegraphics[scale=.18]{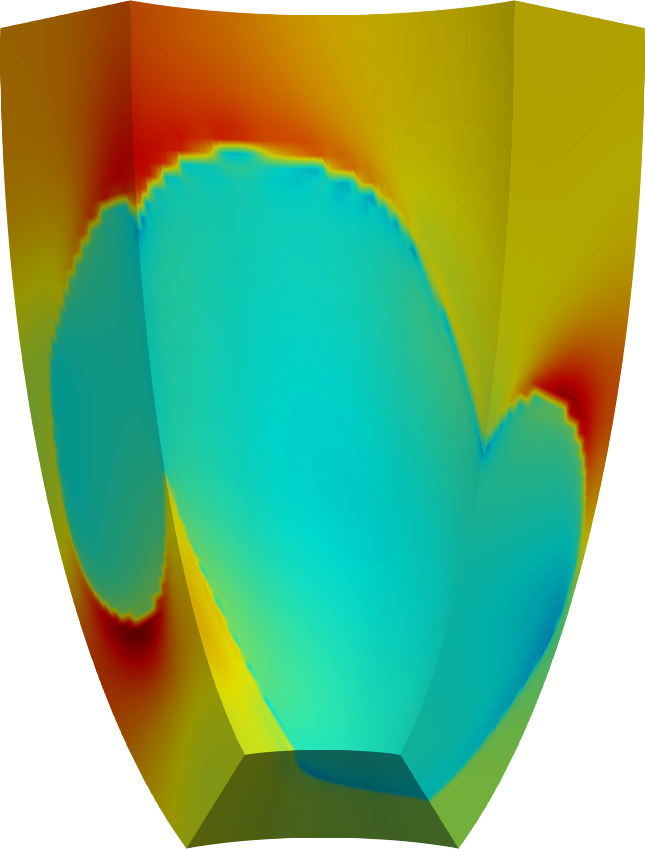} \\	\vspace*{2mm}
		\end{multicols}
		\includegraphics[scale=.50]{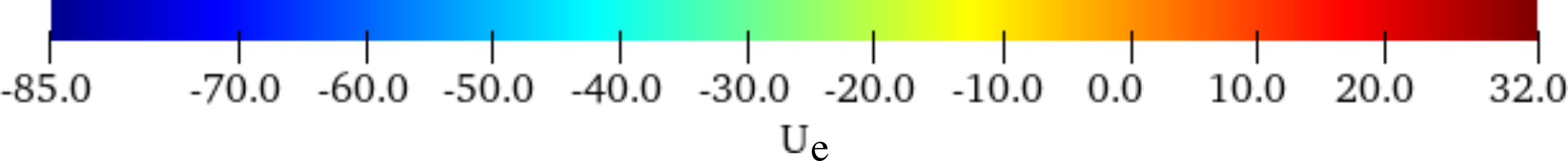}
	\end{center}
	\caption{\emph{Snapshots of the transmembrane $v$ and extracellular $u_e$ potentials at $1$ ms, $20$ ms and $40$ ms.} For each time frame, we report the endocardial view of a portion of the left ventricle, modeled as a truncated ellipsoid. }
	\label{fig: fig_bido_snapshots}
\end{figure}

\begin{figure}[pb]
	\centering
	\begin{subfigure}{\textwidth}
		\begin{tikzpicture}
			\begin{semilogyaxis}
				[width=0.5\textwidth, height=\figHeight, no markers, title, xlabel=Time $(ms)$, ylabel=Nonlinear iterations, tick label style={font=\tiny},line width=1.5pt,legend style={draw=none, fill opacity=0.7, text opacity = 1,row sep=2pt, font=\small, at={(0.98, 0.5,)}, anchor=east}]
				\addplot+[\newtoncolor]  table [x=Time, y=snesIts, col sep=comma, each nth point=200, filter discard warning=false, unbounded coords=discard] {results/heartbeat/HEARTBEAT_NEWTON_32.csv};
				\addplot+[\inewtoncolor]  table [x=Time, y=snesIts, col sep=comma, each nth point=200, filter discard warning=false, unbounded coords=discard ] {results/heartbeat/HEARTBEAT_NEWTON_EW_32.csv};
				\addplot+[\bfgscolor]  table [x=Time, y=snesIts, col sep=comma, each nth point=200, filter discard warning=false, unbounded coords=discard ] {results/heartbeat/HEARTBEAT_QN_PREONLY_32.csv};
				\addplot+[\ibfgscolor]  table [x=Time, y=snesIts, col sep=comma, each nth point=200, filter discard warning=false, unbounded coords=discard ] {results/heartbeat/HEARTBEAT_QN_JACOBIAN_LOW_32.csv};
				\legend {Newton-MG, iNewton, QN preonly, QN jac-low}
			\end{semilogyaxis}
		\end{tikzpicture}
		\hspace{.05\linewidth}
		\begin{tikzpicture}
			\begin{axis}
				[width=0.475\textwidth, height=\figHeight, no markers, title=, xlabel=Time $(ms)$, ylabel=, tick label style={font=\tiny},line width=1.5pt,legend style={draw=none, fill opacity=0.7, text opacity = 1,row sep=2pt, font=\small}, legend pos=north west, legend pos=north east]
				\addplot+[\ngmrescolor]  table [x=Time, y=snesIts, col sep=comma, each nth point=200, filter discard warning=false, unbounded coords=discard ] {results/heartbeat/HEARTBEAT_NGMRES_32.csv};
				\addplot+[\ncgcolor]  table [x=Time, y=snesIts, col sep=comma, each nth point=200, filter discard warning=false, unbounded coords=discard ] {results/heartbeat/HEARTBEAT_NCG_32.csv};
				\legend {NGMRES, NCG-PRP}
			\end{axis}
		\end{tikzpicture}
	\end{subfigure}
	
	\begin{subfigure}{\textwidth}
		\begin{tikzpicture}
			\begin{axis}
				[width=0.98\textwidth, height=\figHeight, no markers, title=, xlabel=Time $(ms)$, ylabel=Solution time (s), tick label style={font=\tiny},line width=1.5pt,legend style={draw=none, fill opacity=0.7, text opacity = 1,row sep=2pt, font=\small}, legend pos=north east, legend columns=3]
				\addplot+[\newtoncolor] table [x=Time, y=SNEStime, col sep=comma, each nth point=200, filter discard warning=false, unbounded coords=discard ] {results/heartbeat/HEARTBEAT_NEWTON_32.csv};
				\addplot+[\inewtoncolor]  table [x=Time, y=SNEStime, col sep=comma, each nth point=200, filter discard warning=false, unbounded coords=discard ] {results/heartbeat/HEARTBEAT_NEWTON_EW_32.csv};
				\addplot+[\bfgscolor]  table [x=Time, y=SNEStime, col sep=comma, each nth point=200, filter discard warning=false, unbounded coords=discard ] {results/heartbeat/HEARTBEAT_QN_PREONLY_32.csv};
				\addplot+[\ibfgscolor]  table [x=Time, y=SNEStime, col sep=comma, each nth point=200, filter discard warning=false, unbounded coords=discard ] {results/heartbeat/HEARTBEAT_QN_JACOBIAN_LOW_32.csv};
				\addplot+[\ngmrescolor]  table [x=Time, y=SNEStime, col sep=comma, each nth point=200, filter discard warning=false, unbounded coords=discard ] {results/heartbeat/HEARTBEAT_NGMRES_32.csv};
				\addplot+[\ncgcolor]  table [x=Time, y=SNEStime, col sep=comma, each nth point=200, filter discard warning=false, unbounded coords=discard ] {results/heartbeat/HEARTBEAT_NCG_32.csv};
				\legend {Newton-MG, iNewton, QN preonly, QN jac-low, NGMRES, NCG-PRP}
			\end{axis}
		\end{tikzpicture}
	\end{subfigure}
	\caption{\emph{Full activation-recovery simulation.} The performance of the different nonlinear solvers is separated into Newton (left) and first order methods (right) as the iteration numbers present different orders of magnitude. On the first row we show the nonlinear iterations incurred by each method, whereas on the second row we show the total solution time of each instant in seconds. The nonlinear CG method is plotted until roughly $200$ ms, after which it exceeded the maximum number of nonlinear iterations allowed, set to 2000.}
	\label{fig:heartbeat}
\end{figure}
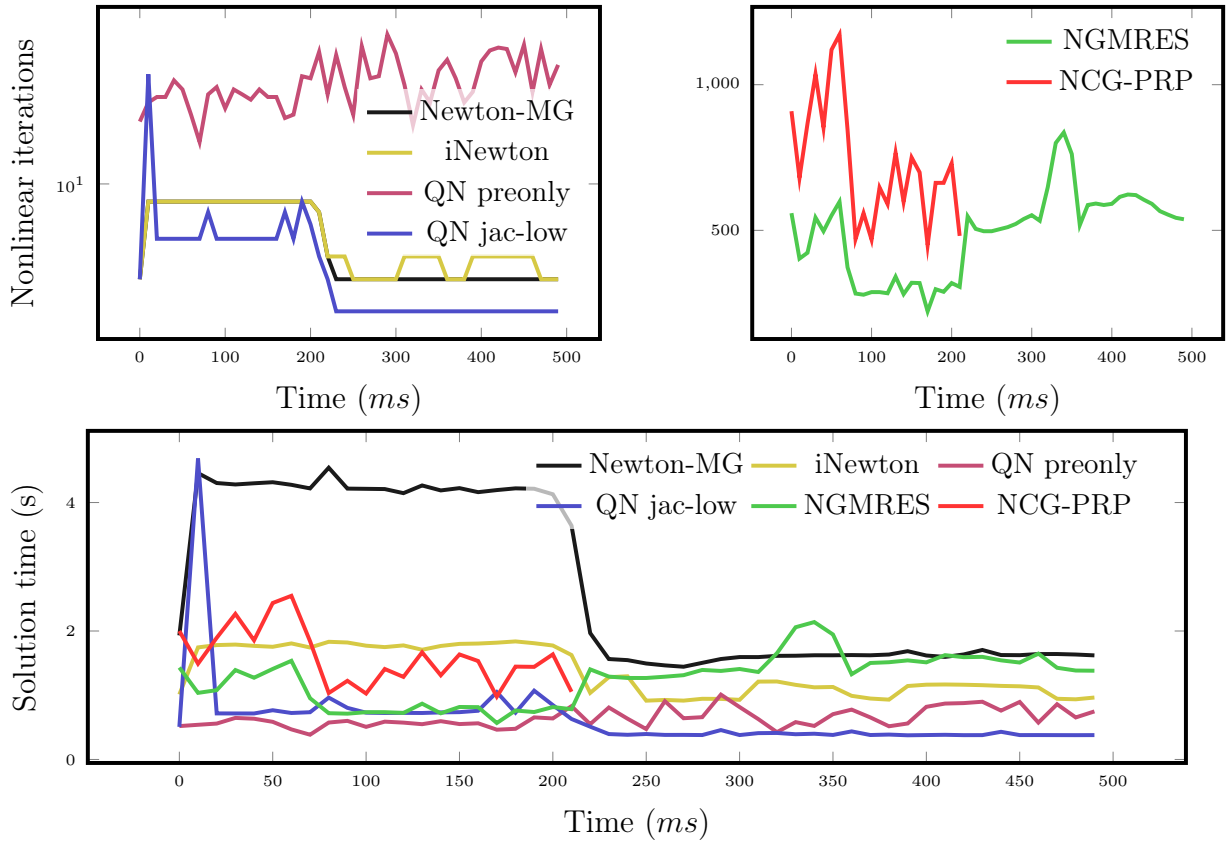

We now compare the performance of the nonlinear solvers during a full activation-recovery phase, by considering a time interval of $[0,500]$ ms (10'000 time steps). The number of processors is fixed to $N_p=12$ and the mesh consists of $ 32\times 32\times 32$ elements, resulting in approximatively $72$ thousand degrees of freedom. Results are shown in Figure \ref{fig:heartbeat}, while Figure \ref{fig: fig_bido_snapshots} shows the transmembrane $v$ and extracellular $u_e$ potentials at different time frames.

We first observe that NCG-PRP is plotted up to roughly $200$ ms because after that it exceeded the maximum number of nonlinear iterations allowed, which we set to 2'000. Besides this method, all Newton methods (inexact, quasi and standard) present a small deterioration during the activation at the beginning of the simulation, and then present a robust behavior. 
This deterioration is instead less appreciated in the first order methods, where NGMRES deteriorates more in the interval $[300, 400]$ ms, where indeed it performs even worse than Newton-MG. Overall, we observe that all methods outperform Newton-MG. 
In particular, the best performance is obtained by the QN preonly method during the first $200$ ms, and by the QN jac-low method after that moment, concluding that quasi-Newton present the best performance. Still, we highlight that during the first $200$ ms, NGMRES outperforms both Newton-MG and inexact-Newton methods, and indeed during almost half of that interval it performs as well as the QN jac-low method. 
All these behaviors can be partially explained by the physical scenario, since the portion of tissue represented in the simulation is almost at rest after 200 ms. Further studies should investigate an adaptive choice of the method parameters depending on the different phases of the electric propagation.

\subsection{Scalability}
We then verify the strong scalability of the proposed methods, by testing our Bidomain solvers with roughly 1.8 million degrees of freedom with up to 128 cores. The results are shown in Table \ref{tab:scalability}, where the speed-up and parallel efficiency are computed respectively as
$$ S_p = \frac{T_1}{T_p} \quad\text{and}\quad E_p = \frac{T_1}{p \,T_p}, $$
where $T_\ell$ denotes the time required by a given solver when using $\ell$ processors. We note that in this case the speed-up is not necessarily a good indicator: as a matter of fact, it is possible to rank the methods according to their CPU time and according to their speed-up, both considering 128 cores -- and the results are not in agreement, as shown in Table \ref{tab:scalability ranking}. 
\begin{table}[pb]
	\caption{\emph{Strong scalability.} Ranking of the nonlinear solvers according to their speed-up and CPU time, when considering 128 cores.}
	\label{tab:scalability ranking}
	\centering
	\begin{tabular}{@{}rcc@{}}
		\toprule rank & Speed-up & CPU time \\ \midrule
		\# 1 & NCG-PRP        & QN preonly \\
		\# 2 & NGMRES     & QN jac-low \\
		\# 3 & Newton-MG  & iNewton \\
		\# 4 & iNewton 	& Newton-MG\\
		\# 5 & QN preonly & NGMRES\\
		\# 6 & QN jac-low & NCG-PRP \\\bottomrule
	\end{tabular}
\end{table}
The results are almost inverted, and in this case the speed-up can be a deceitful indicator. Additionally, this comparison reinforces the superiority of quasi-Newton methods, since they provide the lowest CPU time. We also observe that first order methods present exactly the same number of nonlinear iterations, which can be motivated by the absence of a preconditioner, which in turns allows them to be exactly the same method in parallel, so that their speed-up is that of parallelizing level 1 BLAS operations (operations between vectors). 
Lastly, at first sight it might seem surprising that first order methods are slower than the Newton-MG method, which is in contrast to the results of Table \ref{tab:sensitivity}. This can be explained by the elevated number of degrees of freedom used, which deteriorates the performance of first order methods as shown in Figure \ref{fig:robustness ngmres ncg}.

We highlight that all communication is done through the use of PETSc operations, which are highly optimized for parallel computing. This justifies why we observe an adequate speed-up only up to 4 processors, where there is an average load per CPU of roughly half a million degrees of freedom. Further tuning of the preconditioners to improve this behavior could have been done, but was beyond the scope of this work. The value of these results is mainly comparative among them. 

\begin{table}
	\caption{\emph{Strong scalability and speed-up.} Comparison of the performance between all the methods, in terms of computational times in seconds (CPU), number of nonlinear iterations (NL), speed-up ($S_p$) and parallel efficiency ($E_p$). Fixed problem size of 1.8 million DoFs in the first two time steps, increasing number of cores from 1 to 128.}
	\label{tab:scalability}
	\centering
	\begin{tabular}{@{}rrrrrrrrr@{}}
		\toprule & \multicolumn{4}{c}{Newton-MG} & \multicolumn{4}{c}{iNewton}  \\ \midrule
		Cores & CPU  & NL & $S_p$ & $E_p$ & CPU  & NL & $S_p$ & $E_p$ \\ \midrule
		1    & 485.72 & 3.0   & 1.00  & 1.00 & 307.73 & 4.0   & 1.0   & 1.00 \\
		2    & 239.53 & 3.0   & 2.02  & 1.01 & 163.40 & 4.0   & 1.87  & 0.94 \\
		4    & 127.65 & 3.0   & 3.80  & 0.95 & 87.64  & 4.0   & 3.51  & 0.88 \\
		8    & 77.67  & 3.0   & 6.25  & 0.78 & 45.31  & 3.5   & 6.79  & 0.85 \\
		16   & 55.31  & 3.0   & 8.78  & 0.54 & 34.52  & 4.0   & 8.91  & 0.56 \\
		32   & 55.63  & 3.0   & 8.73  & 0.27 & 29.00  & 4.0   & 10.61 & 0.33 \\
		64   & 30.88  & 3.0   & 15.73 & 0.25 & 19.29  & 4.0   & 15.95 & 0.25 \\
		128  & 17.03  & 3.0   & 28.52 & 0.22 & 11.03  & 4.0   & 27.90 & 0.22 \\ \midrule 
		\toprule & \multicolumn{4}{c}{QN jac-low} & \multicolumn{4}{c}{QN preonly} \\ \midrule
		Cores & CPU  & NL & $S_p$ & $E_p$ & CPU  & NL & $S_p$ & $E_p$ \\ \midrule
		1    & 485.72 & 3.0   & 1.00  & 1.00  & 122.65 & 32.0  & 1.00  & 1.00 \\
		2    & 239.53 & 3.0   & 2.02  & 1.01  & 97.18  & 55.5  & 1.26  & 0.63 \\
		4    & 127.65 & 3.0   & 3.80  & 0.95  & 37.34  & 34.0  & 3.28  & 0.82 \\
		8    & 77.67  & 3.0   & 6.25  & 0.78  & 21.02  & 32.0  & 5.83  & 0.73 \\
		16   & 55.31  & 3.0   & 8.78  & 0.54  & 17.42  & 35.0  & 7.04  & 0.44 \\
		32   & 55.63  & 3.0   & 8.73  & 0.27  & 13.50  & 32.5  & 9.09  & 0.28 \\
		64   & 30.88  & 3.0   & 15.73 & 0.25  & 9.57   & 31.0  & 12.82 & 0.20 \\
		128  & 17.03  & 3.0   & 28.52 & 0.22  & 5.76   & 32.5  & 21.29 & 0.17 \\ \midrule 
		\toprule & \multicolumn{4}{c}{NGMRES} & \multicolumn{4}{c}{NCG-PRP} \\ \midrule
		Cores	& CPU  & NL & $S_p$ & $E_p$ & CPU & NL & $S_p$ & $E_p$ \\ \midrule
		1   & 658.70 & 524.0 & 1.0   & 1.00 & 958.08 & 956.0 & 1.00  & 1.00\\
		2   & 320.29 & 524.0 & 2.06  & 1.03 & 507.45 & 956.0 & 1.89  & 0.94 \\
		4   & 157.48 & 524.0 & 4.18  & 1.05 & 232.67 & 956.0 & 4.12  & 1.03\\
		8   & 92.59  & 524.0 & 7.11  & 0.89 & 137.00 & 956.0 & 6.99  & 0.87\\
		16  & 70.95  & 524.0 & 9.28  & 0.58 & 103.23 & 956.0 & 9.28  & 0.58\\
		32  & 70.60  & 524.0 & 9.33  & 0.29 & 93.39  & 956.0 & 10.26 & 0.32\\
		64  & 31.86  & 524.0 & 20.67 & 0.32 & 54.74  & 956.0 & 17.50 & 0.27 \\
		128 & 19.32  & 524.0 & 34.09 & 0.27 & 27.89  & 956.0 & 34.35 & 0.27\\ \bottomrule
	\end{tabular}
\end{table}

\subsection{Convergence}
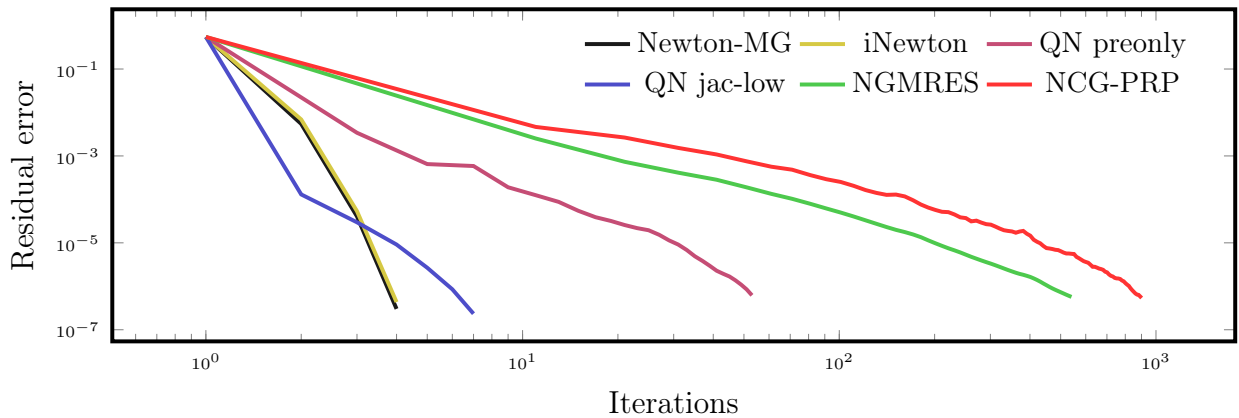
\begin{figure}
	\begin{tikzpicture}
		\begin{loglogaxis}
			[width=\linewidth, height=\figHeight, no markers, xlabel=Iterations, ylabel=Residual error, tick label style={font=\tiny},line width=1.5pt,legend style={draw=none, fill opacity=0.7, text opacity = 1,row sep=2pt,font=\small}, legend pos=north east, legend columns=3, domain=10:100]
			\addplot+[\newtoncolor]  table [x=its, y=Newton, col sep=comma] {convergence.csv};
			\addplot+[\inewtoncolor]  table [x=its, y=iNewton, col sep=comma] {convergence.csv};
			\addplot+[\bfgscolor]  table [x=its, y=QN, col sep=comma, each nth point=2] {convergence.csv};
			\addplot+[\ibfgscolor]  table [x=its, y=iQN, col sep=comma] {convergence.csv};
			\addplot+[\ngmrescolor]  table [x=its, y=NGMRES, col sep=comma, each nth point=10] {convergence.csv};
			\addplot+[\ncgcolor]  table [x=its, y=NCG, col sep=comma, each nth point=10] {convergence.csv};
			\legend{Newton-MG, iNewton, QN preonly, QN jac-low, NGMRES, NCG-PRP};
		\end{loglogaxis}
	\end{tikzpicture}
	\caption{\emph{Convergence test.} Evolution of the residual errors at the first time step for all solvers considered.}
	\label{fig:errors}
\end{figure}

In this section, we analyze the evolution of the residual errors for each of the Bidomain solvers considered, collecting the results in Figure \ref{fig:errors}. We make two observations. (i) After the first few iterations, the curve slopes roughly represent the rates of convergence of the considered methods, meaning that Newton-MG and inexact-Newton have the steepest slopes (faster convergence), accordingly to the theoretical expectations; first order methods present less steep slopes, but similar among NGMRES and NCG-PRP, while quasi-Newton methods have intermediate slopes. (ii) Quasi-Newton methods present a very fast decay of the error in the first iterations, even more than the standard Newton: this fact could inspire the usage of hybrid nonlinear solvers, where different methods are used at different phases of the iterations. This has not been considered in literature and could be an interesting road to pursue in future research. 

\rev{
\subsection{Comparison with IMEX discretization}
The proposed methods all depend on an implicit treatment of the nonlinear term $I_\text{ion}$, which as we have noted is useful for high order methods, but it is not what is commonly used in practice. Indeed, most electrophysiology solvers \cite{lifex-ep, opencarp} use (i) an explicit treatment of the nonlinear term and (ii) a pseudo-Bidomain time discretization, where the parabolic-elliptic formulation is used, and the extracellular potential $u_e$ is computed in a staggered fashion. For this reason, we would like to give evidence that, in a practical setting, it can be better to use an implicit instead of explicit treatment of the nonlinear term in the parabolic equation. To this end, we consider the simplified scenario of the monodomain model with the FitzHugh-Nagumo model to compare both implicit and explicit treatments of the nonlinearity in an idealized left ventricle, displayed in Figure~\ref{fig:heart}. The geometry is non-structured, and the solvers were implemented using Firedrake \cite{firedrake} in an in-house library that is currently under development.

\begin{figure}[ht!]
    \centering
    \includegraphics[width=0.4\textwidth]{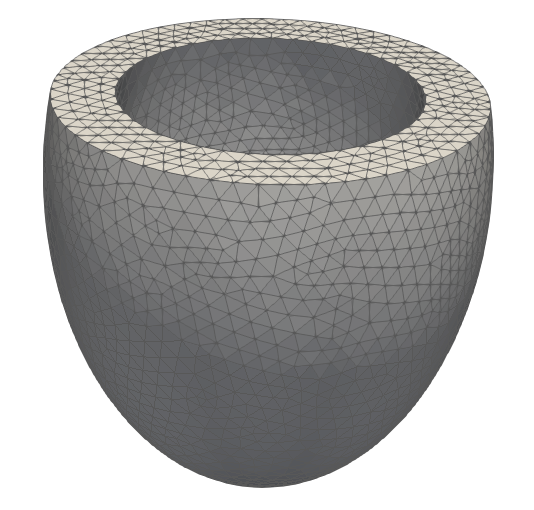}
    \caption{Idealized ventricle geometry used to compare IMEX and implicit discretizations.}
    \label{fig:heart}
\end{figure}

To make the comparison fair, we have first computed the factor $N$ such that the IMEX discretization is as accurate as the implicit one. To do this, we fixed the time step size of the implicit method at $\tau=0.02\,ms$. Then we computed an exact-in-time solution for a timestep given by $\tau/1000$, and we computed the Bochner norm of the error during the first $50$ timesteps of simulation, according to 
    $$ \|\vec x\|^2_{L^2(0,T; H^1(\Omega))} \approx \tau\sum_n \|\vec x^n\|^2_{H^1(\Omega)}.$$
Doing so, we obtained that the factor is given by $N=2$, i.e. the accuracy of the IMEX discretization with $\tau/2$ is roughly the same as the implicit one using $\tau$. The described problem was then simulated for $20\,ms$ in a single CPU core, from which we obtained the CPU times per timestep shown in Table~\ref{table:imex}. To make this comparison as fair as possible, we used the best available methods for each formulation, meaning that we fixed the preconditioner to be an algebraic multigrid. Then, we solved the IMEX linear system using the Conjugate Gradient method, and the implicit one with the BFGS preonly method. We highlight that, even though the IMEX scheme yields lower CPU times per timestep, if we fix the accuracy of both methods, then the implicit solver using the methods developed in this work yield the same accuracy in a $55\%$ of the time. In the future we will expand this approach to different and more complex models to define what is the best solver for each model.

\begin{table}[ht!]
    \centering
    \begin{tabular}{c | c c c}
    \toprule  & CPU time per PDE solve & CPU time per timestep & Total CPU time for $20\,ms$ \\\midrule
    IMEX      & $19.73\,ms$        & $24.55\,ms$     & $24.55\,s$ \\
    Implicit  & $22.33\,ms$        & $27.11\,ms$     & $13.56\,s$\\\bottomrule
    \end{tabular}
    \caption{Performance comparison between IMEX and implicit time schemes. In the first column, we show the time it takes to solve the corresponding PDE once. In the second column, we show the time it takes to solve a total time of $\tau$, considering both ODE and PDE solutions. The third column shows the total time it takes to simulate $20\,ms$ in our implementation.}
    \label{table:imex}
\end{table}
}

\section{Conclusions} \label{sec: conclusions}
The computation of a variational principle for a decoupled implicit time discretization of the Bidomain equations allowed us to analytically prove the convergence of two classes of nonlinear solvers, namely quasi-Newton methods and the nonlinear CG Fletcher-Reeves descent method. For these methods, and also the nonlinear GMRES and other variants of NCG, we have performed a thorough numerical study to verify their robustness and scalability, and all methods but the NCG present very satisfactory results. Besides this, quasi-Newton methods present the best overall performance. Both variants (QN preonly and QN jac-low) are comparable, and a significant improvement (of over 60\% in many cases) should be expected when changing the standard Newton-MG solver. A simpler modification would be to consider an inexact-Newton method, which requires minimal modifications to a working Newton code and has shown to be not only faster, but in many cases also more robust, than the standard Newton method \cite{barnafi2022alternative}. \rev{We have further compared the implicit and IMEX discretizations by using comparable accuracy, which allows us to show that our approach can be superior to the typically used IMEX formulations, potentially halving the overall CPU time and saving communication in the form of vector updates. Still, the scope of this work was not that of providing a detailed comparison of these two schemes, so we expect to extend this investigation in future studies.} 

We remark that we dealt with order one both for time and space discretizations, thus we should expect that a variation of these choices would affect also the behaviour of all the above mentioned nonlinear solvers.

The success of nonlinear GMRES, a first order method, opens the possibility of using a matrix-free approach for electrophysiology, which could be computed exclusively in a GPU. We believe this to be fundamental for future studies of hybrid computing in electromechanics.

\appendix
\section{PETSc instructions}\label{appendix:petsc}
In this appendix, we provide the PETSc command line instructions used for calling each of the methods under consideration.
\begin{itemize}
	\item {\bf Newton-MG }
	\begin{lstlisting}[language=bash, frame=single, caption=PETSc commands to use Newton-MG.]
		-snes_type newtonls
		-ksp_type  cg
		-pc_type   gamg
		-snes_atol  1e-12
		-snes_rtol  1e-6
		-snes_stol  0.0
		-snes_max_it  2000
		-ksp_constant_null_space
	\end{lstlisting}
	\item {\bf Inexact Newton-MG}
	\begin{lstlisting}[language=bash, frame=single, caption=PETSc commands to use inexact Newton-Krylov.]
		-snes_type newtonls
		-ksp_type  cg
		-pc_type   gamg
		-snes_atol  1e-12
		-snes_rtol  1e-6
		-snes_stol  0.0
		-snes_max_it  2000
		-snes_ksp_ew
		-snes_ksp_ew_rtol 1e-1 # Tuned parameter
		-ksp_constant_null_space
		-ksp_atol 0.0
	\end{lstlisting}
	\item {\bf QN jac-low}
	\begin{lstlisting}[escapeinside=``,language=bash, frame=single, caption=PETSc commands to use inexact-BFGS.]
		-snes_type qn
		-ksp_type  cg
		-pc_type   gamg
		-snes_atol  1e-12
		-snes_rtol  1e-6
		-snes_stol  0.0
		-snes_max_it  2000
		-snes_qn_type lbfgs
		-snes_qn_m 10 # Tuned parameter
		-snes_lag_jacobian 9999
		-snes_lag_preconditioner 9999 
		-snes_qn_restart_type none
		-ksp_constant_null_space
		-ksp_norm_type none
		-ksp_max_it 10
	\end{lstlisting}
	\item {\bf QN-preonly}
	\begin{lstlisting}[language=bash, frame=single, caption=PETSc commands to use BFGS.]
		-snes_type qn
		-ksp_type  preonly
		-pc_type   gamg
		-snes_atol  1e-12
		-snes_rtol  1e-6
		-snes_stol  0.0
		-snes_qn_type lbfgs
		-snes_qn_m 10 # Tuned parameter
		-snes_qn_scale_type jacobian
		-snes_lag_jacobian 9999
		-snes_lag_preconditioner 9999 
		-snes_qn_restart_type none
		-ksp_constant_null_space
	\end{lstlisting}
	\item {\bf Nonlinear GMRES}
	\begin{lstlisting}[language=bash, frame=single, caption=PETSc commands to use NGMRES.]
		-snes_type ngmres
		-snes_atol  1e-12
		-snes_rtol  1e-6
		-snes_stol  0.0
		-snes_max_it  2000
		-snes_ngmres_restart_type none
		-snes_ngmres_m 10 # Tuned parameter
	\end{lstlisting}
	\item {\bf Nonlinear CG-PRP}
	\begin{lstlisting}[language=bash, frame=single, caption=PETSc commands to use NCG-PRP.]
		-snes_type ncg
		-snes_atol  1e-12
		-snes_rtol  1e-6
		-snes_stol  0.0
		-snes_max_it  2000
		-snes_ncg_type prp # Tuned parameter
	\end{lstlisting}
	
\end{itemize}

\section*{Acknowledgements}
N.A. Barnafi, N.M.M. Huynh and L.F. Pavarino have been supported by grants of MIUR (PRIN 2017AXL54F$\_$002) and INdAM--GNCS. 
N.A. Barnafi and S. Scacchi have been supported by grants of MIUR (PRIN 2017AXL54F$\_$003) and INdAM-GNCS. 
N.M.M. Huynh, L.F. Pavarino and S. Scacchi have been supported by the European High-Performance Computing Joint Undertaking EuroHPC under grant agreement No 955495 (MICROCARD) co-funded by the Horizon 2020 programme of the European Union (EU), and the Italian ministry of economic development.
N.A. Barnafi was supported by the ANID Grant \emph{FONDECYT de Postdoctorado N° 3230326} and by Centro de Modelamiento Matemático, Proyecto Basal FB210005. 
The Authors are also grateful to the University of Pavia for the usage of the cluster EOS.


\end{document}

%% file: mycustom.tex
\usepackage{algorithm}
\usepackage{algpseudocode}
\usepackage{amsmath,amssymb,amsfonts}
\usepackage{booktabs}
\usepackage{hyperref} 
\usepackage{mathtools}
\usepackage{multicol}
\usepackage{multirow}
\usepackage{pgfplots}
	\usepgfplotslibrary{colorbrewer}
	\pgfplotsset{cycle list/Set2}
\usepackage{subcaption}
\usepackage[T1]{fontenc}
\usepackage{fullpage}

\pgfplotsset{compat=1.18}

\usepackage{isomath}
\renewcommand{\vec}{\vectorsym}
\newcommand{\mat}{\matrixsym}
\newcommand{\ten}{\tensorsym}

\DeclareMathOperator{\dive}{div}
\DeclareMathOperator{\grad}{\nabla}
\DeclareMathOperator{\argmin}{argmin}

\newcommand{\xx}{\mathbf x}
\newcommand{\aal}{\mathbf a_l}
\newcommand{\aat}{\mathbf a_t}
\newcommand{\aan}{\mathbf a_n}

\newtheorem{theorem}{Theorem}
\newtheorem{remark}{Remark}

\newcommand{\Hone}{V_0} 
\newcommand{\R}{\mathbb R}
\newcommand{\normL}[1]{\|#1\|_{L^2(\Omega)}}

\newcommand{\seminormH}[1]{|#1|_{H^1(\Omega)}}

\usepackage{listings}
\usepackage{xcolor}
\definecolor{codegreen}{rgb}{0,0.6,0}
\definecolor{codegray}{rgb}{0.5,0.5,0.5}
\definecolor{codepurple}{rgb}{0.58,0,0.82}
\definecolor{backcolour}{rgb}{0.95,0.95,0.92}
\lstdefinestyle{mystyle}{
	backgroundcolor=\color{backcolour}, commentstyle=\color{codegreen},
	keywordstyle=\color{magenta},
	numberstyle=\tiny\color{codegray},
	stringstyle=\color{codepurple},
	basicstyle=\ttfamily\footnotesize,
	breakatwhitespace=false,         
	captionpos=b,                    
	keepspaces=true,                 
	numbers=none,                    
	numbersep=5pt,                  
	showspaces=false,                
	showstringspaces=false,
	showtabs=false,                  
	tabsize=2,
	framerule=1.5pt,
	rulecolor=\color{red!60!black},
	float
}
\newcommand{\newtoncolor}{black!90!white}
\newcommand{\inewtoncolor}{yellow!80!black}
\newcommand{\bfgscolor}{purple!90!black!70!white}
\newcommand{\ibfgscolor}{blue!70!black!70!white}
\newcommand{\ngmrescolor}{green!70!black!70!white}
\newcommand{\ncgcolor}{red!80!white}

\newcommand{\rev}[1]{{\color{red}#1}}

%% file: main.bbl
\begin{thebibliography}{99}
	
    \bibitem{lifex-ep}
    P.C. Africa, R. Piersanti, F. Regazzoni, M. Bucelli, M. Salvador, M. Fedele, S. Pagani, L. Dede’, A. Quarteroni, 
    lifex-ep: a robust and efficient software for cardiac electrophysiology simulations,
    \emph{BMC bioinformatics}, \textbf{24} (2023), 389.

	\bibitem {al-baali1985descent}
	M. Al-Baali, 
	Descent property and global convergence of the Fletcher-Reeves method with inexact line search, 
	\emph{I.M.A. J. Numer. Analysis}, \textbf{5} (1985), pp. 121–124.
	
	\bibitem{balay2019petsc}
	S. Balay et al.,
	{PETS}c users manual,
	(2019).
	
	\bibitem{barnafi2021mathematical}
	N. Barnafi, P. Zunino, L. Ded{\'e} and A.M. Quarteroni,
	Mathematical analysis and numerical approximation of a general linearized poro-hyperelastic model,
	\emph{Comp. Math. Appl.}, \textbf{91} (2021), pp. 202--228.
	
	\bibitem{barnafi2022alternative}
	N. Barnafi, N.M.M. Huynh, L.F. Pavarino and S.Scacchi,
	Alternative parallel nonlinear solvers in cardiac modeling,
	\emph{IFAC-PapersOnLine} (2022), 50.20 (2022), pp. 187--192.
	
	\bibitem{barnafi2022parallel}
	N. Barnafi, L.F. Pavarino and S.Scacchi,
	Parallel inexact Newton-Krylov and quasi-Newton solvers for nonlinear elasticity,
	\emph{Comp. Meth. Appl. Mech. Engrg.}, \textbf{400} (2022), pp. 115557.
	
	\bibitem{bjornsson2020digital}
	B. Bj{\"o}rnsson et al., 
	Digital twins to personalize medicine,
	\emph{Genome medicine}, \textbf{12}(1) (2020), pp.1--4. 
	
	\bibitem{bourgault2003simulation}
	Y. Bourgault, M. Ethier and V.G. LeBlanc, 
	Simulation of electrophysiological waves with an unstructured finite element method,
	\emph{ESAIM: Math. Model. Num. Anal.}, \textbf{37}(4) (2003), pp. 649--661.
	
	\bibitem{brune2015composing}
	P.R. Brune, M.G. Knepley, B.F. Smith and X. Tu, 
	Composing scalable nonlinear algebraic solvers,
	\emph{SIAM Review}, \textbf{57}(4) (2015), pp. 535--565.
	
	\bibitem{chen2019splitting}
	H. Chen, X. Li and Y. Wang,
	A two-parameter modified splitting preconditioner for the Bidomain equations,
	\emph{Calcolo}, \textbf{56}(2) (2019), 21.
	
	\bibitem{cornejo2015analisis}
	N.A. Cornejo Fuenzalida,
	Análisis variacional de las ecuaciones de FitzHugh-Nagumo en electrofisiologıa cardıaca,
	\emph{Diss. Pontificia Universidad Católica de Chile} (2015).
	
	\bibitem{franzone2006computational}
	P. Colli Franzone, L.F. Pavarino and G. Savar{\'e},
	Computational electrocardiology: mathematical and numerical modeling,
	\emph{Complex systems in Biomedicine}, Springer (2006), pp. 187--241.
	
	\bibitem{franzone2014mathematical}
	P. Colli Franzone, L.F. Pavarino and S. Scacchi,
	Mathematical cardiac electrophysiology, 
	\emph{Springer},\textbf{13} (2014).
	
	\bibitem{colli2018numerical}
	P. Colli Franzone, L.F. Pavarino and S. Scacchi,
	A numerical study of scalable cardiac electro-mechanical solvers on {HPC} architectures,
	\emph{Front. Physiol.}, \textbf{9} (2018), pp. 268.
	
	\bibitem{dacorogna2007direct}
	B. Dacorogna,
	Direct methods in the calculus of variations,
	\emph{Springer Science \& Business Media}, \textbf{78} (2007).
	
	\bibitem{dede2019computational}
	L. Ded{\'e}, F. Menghini and A.M. Quarteroni,
	Computational fluid dynamics of blood flow in an idealized left human heart,
	\emph{Int. J. Num. Meth. Biomed. Engrg.} (2019), pp. e3287.
	
	\bibitem{di2021computational}
	S. Di Gregorio, M. Fedele, G. Pontone, A.F. Corno, P. Zunino, C. Vergara and A.M. Quarteroni,
	A computational model applied to myocardial perfusion in the human heart: from large coronaries to microvasculature,
	\emph{J. Comput. Physics}, \textbf{424} (2021), pp. 109836.
	
	\bibitem{eisenstat1994globally}
	S.C. Eisenstat and H.F. Walker,
	Globally convergent inexact Newton methods,
	\emph{SIAM J. Opt.}, \textbf{4}(2) (1994), pp. 393--422.
	
	\bibitem{eisenstat1996choosing}
	S.C. Eisenstat and H.F. Walker,
	Choosing the forcing terms in an inexact Newton method,
	\emph{SIAM J. Sci. Comput.}, \textbf{17}(1) (1996), pp. 16--32.
	
	\bibitem{fletcher1964function}
	R. Fletcher and C.M. Reeves, 
	Function minimization by conjugate gradients, 
	\emph{Computer Journal}, \textbf{7} (1964), pp. 149–154.
	
	\bibitem{gelfand2000calculus}
	I.M. Gelfand and R.A. Silverman,
	Calculus of variations,
	\emph{Courier Corporation} (2000).

    \bibitem{griewank1987}
    A. Griewank, 
    The local convergence of Broyden-like methods on Lipschitzian problems in Hilbert spaces, 
    \emph{SIAM J. Num. An.}, \textbf{24}-3 (1987), pp. 684-705.
	
	\bibitem{huynh2021scalable}
	N.M.M. Huynh, L.F. Pavarino and S.Scacchi,
	Scalable Newton-Krylov-BDDC and FETI-DP deluxe solvers for decoupled cardiac reaction-diffusion models,
	\emph{14th WCCM-ECCOMAS Congress 2020}, (2021) \textbf{400}.
	
	\bibitem{huynh2021newton}
	N.M.M. Huynh,
	Newton-Krylov-BDDC deluxe solvers for non-symmetric fully implicit time discretizations of the Bidomain model,
	\emph{Numerische Mathematik}, 152(4) (2022), pp. 841--879.
	
	\bibitem{huynh2022parallel}
	N.M.M. Huynh, L.F. Pavarino and S. Scacchi,
	Parallel Newton-Krylov-BDDC and FETI-DP deluxe solvers for implicit time discretizations of the cardiac Bidomain equations,
	\emph{SIAM J. Sci. Comput.}, \textbf{44}-2 (2022), pp. B224--B249.
	
	\bibitem{huynh2022scalable}
	N.M.M. Huynh, L.F. Pavarino and S. Scacchi,
	Scalable and robust dual-primal Newton-Krylov deluxe solvers for cardiac electrophysiology with biophysical ionic models,
	\emph{Vietnam J. Math.}50(4) (2022), pp. 1029--1052.
	
	\bibitem{hurtado2014gradient}
	D. Hurtado and D. Henao,
	Gradient flows and variational principles for cardiac electrophysiology: toward efficient and robust numerical simulations of the electrical activity of the heart,
	\emph{Comp. Meth. Appl. Mech. Engrg.}, \textbf{273} (2014), pp.238--254.

	\bibitem{kunisch2013optimal}
	K. Kunisch and M. Wagner
	Optimal control of the bidomain system (ii): Uniqueness and regularity theorems for weak solutions
	\emph{Ann. Mat. Pura Appl.	}, \textbf{192} (2013), pp.951--986.

	\bibitem{legrice1995laminar}
	I.J. LeGrice, B.H. Smaill, L.Z. Chai, S.G. Edgar, J.B. Gavin and P.J. Hunter,
	Laminar structure of the heart: ventricular myocyte arrangement and connective tissue architecture in the dog,
	\emph{Amer. J. Physiol.-Heart Circ.Physiol.}, \textbf{269}-2 (1995), pp. H571--H582.
	
	\bibitem{linge2005solving}
	S. Linge, G. Lines and  J. Sundnes, 
	Solving the heart mechanics equations with Newton and quasi Newton methods--a comparison,
	\emph{Comp. Meth. Biomech. Biomed. Engrg.}, \textbf{8}-1 (2005), pp. 31--38.
	
	\bibitem{liu2017quasi}
	T. Liu, S. Bouaziz and L. Kayan,
	Quasi-Newton methods for real-time simulation of hyperelastic materials,
	\emph{ACM Transactions on Graphics}, \textbf{36}-3 (2017), pp. 1--16

	\bibitem{marsh2012secrets}
	M.E. Marsh, S.T. Ziaratgahi and R.J. Spiteri,
	The secrets to the success of the Rush--Larsen method and its generalizations
	\emph{IEEE Trans. Biomed. Eng.	}, \textbf{59}-9 (2012), pp. 2506--2515.

	\bibitem{munteanu2009decoupled}
	M. Munteanu and L.F. Pavarino,
	Decoupled Schwarz algorithms for implicit discretizations of nonlinear Monodomain and Bidomain systems,
	\emph{Math. Models Methods Appl. Sci.}, \textbf{19}-7 (2009), pp. 1065--1097.
	
	\bibitem{munteanu2009sisc}
	M. Munteanu, L.F. Pavarino and S. Scacchi,
	A scalable Newton–Krylov–Schwarz method for the Bidomain reaction-diffusion system,
	\emph{SIAM J. Sci. Comp.}, \textbf{31}(5) (2009), pp. 3861--3883.
	
	\bibitem{murillo2004fully}
	M. Murillo and X-C. Cai,
	A fully implicit parallel algorithm for simulating the non-linear electrical activity of the heart,
	\emph{Numer. Linear Algebra Appl.}, \textbf{11} (2004), pp. 261--277.

	\bibitem{nagaiah2011numerical}
	C. Nagaiah, K. Junisch and G. Plank
	Numerical solution for optimal control of the reaction-diffusion equations in cardiac electrophysiology
	\emph{Comput. Optim. Appl.}, \textbf{49} (2011), pp. 149--178.

	\bibitem{pennacchio2005multiscale}
	M. Pennacchio, G. Savar{\'e} and P. Colli Franzone,
	Multiscale modeling for the bioelectric activity of the heart,
	\emph{SIAM J. Math. An.}, \textbf{37}(4) (2005), pp. 1333--1370.
	
	\bibitem{pennacchio2011fast}
	M. Pennacchio and V. Simoncini,
	Fast structured amg preconditioning for the Bidomain model in electrocardiology,
	\emph{SIAM J. Sci. Comput.}, \textbf{33}(2) (2011), pp. 721--745.
	
	\bibitem{piersanti20213d}
	R. Piersanti, F. Regazzoni, M. Salvador, A.F. Corno, L. Ded{\'e}, C. Vergara and A.M. Quarteroni,
	3D-0D closed-loop model for the simulation of cardiac biventricular electromechanics,
	\emph{Comp. Meth. Appl. Mech. Engrg.}, \textbf{391} (2021), pp. 114607.

    \bibitem{opencarp}
    G. Plank, A. Loewe, A. Neic, C. Augustin, Y.-L. Huang, M. Gsell, E. Karabelas, M. Nothstein, J. Sánchez, A. Prassl, G. Seemann, and E.Vigmond,
    The {openCARP} Simulation Environment for Cardiac Electrophysiology,
    \emph{Comp. Meth. Pr. Bio.}, \textbf{208} (2021), pp.106223.
	
	\bibitem{quarteroni2017integrated}
	A.M. Quarteroni, T. Lassila, S. Rossi and R. Ruiz-Baier,
	Integrated {H}eart—{C}oupling multiscale and multiphysics models for the simulation of the cardiac function,
	\emph{Comp. Meth. Appl. Mech. Engrg.}, \textbf{314} (2017), pp. 345--407.

    \bibitem{firedrake}
    F. Rathgeber, D.A. Ham, L. Mitchell, M. Lange, F. Luporini, A.T.T McRae, G.-T. Bercea, G.R. Markall, and P.H.J. Kelly,
    Firedrake: automating the finite element method by composing abstractions,
    \emph{ACM TOMS}, \textbf{43} (2016), 3 pp. 1--27.

    \bibitem{sachs1986broyden}
    E.W. Sachs,
    Broyden's method in Hilbert space, 
    \emph{Mathematical Programming}, \textbf{35} (1986), pp. 71--82.
	
	\bibitem{scacchi2011multilevel}
	S. Scacchi,
	A multilevel hybrid Newton-Krylov-Schwarz method for the Bidomain model of electrocardiology,
	\emph{Comp. Meth. Appl. Mech. Engrg.}, \textbf{200}(5--8) (2011), pp. 717--725.
	
	\bibitem{smith2004multiscale}
	N.P. Smith, D.P. Nickerson, E.J. Crampin and P.J. Hunter,
	Multiscale computational modelling of the heart,
	\emph{Acta Numerica}, \textbf{13} (2004), pp. 371--431.
	
	\bibitem{sundnes2005operator}
	J. Sundnes, G. Lines and A. Tveito,
	An operator splitting method for solving the Bidomain equations coupled to a volume conductor model for the torso,
	\emph{Math. Biosciences}, \textbf{194}(2), pp. 233--248.
	
	\bibitem{ten2004model}
	K.H.W.J. Ten Tusscher, D. Noble, P-J. Noble and A.V. Panfilov,
	A model for human ventricular tissue,
	\emph{Amer. J. Physiol.-Heart Circ. Physiol.}, \textbf{286}-4 (2004), pp. H1573--H1589.
	
	\bibitem{veneroni2009reaction}
	M. Veneroni,
	Reaction--diffusion systems for the macroscopic Bidomain model of the cardiac electric field,
	\emph{Nonlinear Anal. Real World Appl.}, \textbf{10}-2 (2009), pp. 849--868.
	
	\bibitem{washio1997krylov}
	T. Washio and C.W. Oosterlee, 
	Krylov subspace acceleration for nonlinear multigrid schemes,
	\emph{Elec. Trans. Num. Anal.}, \textbf{6} (1997), pp. 271--290.

    \bibitem{weiser2005asymptotic}
    M. Weiser, A. Schiela and P. Deuflhard,
    Asymptotic mesh independence of Newton's method revisited,
    \emph{SIAM J. Num. An.}, \textbf{42}-5 (2005), pp. 1830--1845.
	
	\bibitem{wright1999numerical}
	S. Wright and J. Nocedal,
	Numerical optimization,
	\emph{Springer Science}, \textbf{35} (1999).
	
	\bibitem{zampini2014dual}
	S. Zampini,
	Dual-primal methods for the cardiac Bidomain model,
	\emph{Math. Models Methods Appl. Sci.}, \textbf{24}-4 (2014), pp. 667--696.
	
\end{thebibliography}
